\documentclass{siamart251216}
\usepackage{amssymb}
\usepackage{microtype}
\newsiamremark{remark}{Remark}
\crefname{remark}{Remark}{Remarks}

\newcommand{\RK}{Runge--Kutta}
\newcommand{\RKN}{Runge--Kutta--Nystr\"om}
\newcommand{\Firedrake}{\lstinline{Firedrake}}
\newcommand{\Irksome}{\lstinline{Irksome}}

\usepackage{xcolor}
\usepackage{listings}
\definecolor{DarkBlue}{rgb}{0.00,0.00,0.55}
\definecolor{DarkRed}{rgb}{0.55,0.00,0.00}
\definecolor{DarkGreen}{rgb}{0.00,0.55,0.00}
\definecolor{Bittersweet}{rgb}{1.0, 0.44, 0.37}
\definecolor{Purple}{rgb}{0.5, 0.0, 0.5}
\definecolor{seaborngreen}{rgb}{0.3333333333333333, 0.6588235294117647, 0.40784313725490196}  % #55A868
\definecolor{seaborncyan}{rgb}{0.39215686274509803, 0.7098039215686275, 0.803921568627451}  % #64B5CD
\definecolor{seabornblue}{rgb}{0.2980392156862745, 0.4470588235294118, 0.6901960784313725}  % #4C72B0
\definecolor{seabornpurple}{rgb}{0.5058823529411764, 0.4470588235294118, 0.6980392156862745}  % #8172B2
\definecolor{seabornred}{rgb}{0.7686274509803922, 0.3058823529411765, 0.3215686274509804}  % #C44E52
\definecolor{seabornorange}{rgb}{0.958, 0.476, 0.206}  % #F47935
\definecolor{seabornsand}{rgb}{0.8, 0.7254901960784313, 0.4549019607843137}  % #CCB974

\usepackage{pgfplots}
\usepackage{pgfplotstable}
\usepackage[nameinlink]{cleveref}
\usepackage{tikz}
\usepackage{verbatim}
\usepackage{subcaption}
\usepackage{xspace}
\usepackage{todonotes}
\newif\ifdraft
\draftfalse
\ifdraft
  \newcommand{\pef}[1]{\todo[color=blue!30,inline]{{\bf PEF:} #1}}
  \newcommand{\bda}[1]{\todo[color=orange!30,inline]{{\bf BDA:} #1}}
  \newcommand{\rck}[1]{\todo[color=green!30,inline]{{\bf RCK:} #1}}
  \newcommand{\spm}[1]{\todo[color=purple!30,inline]{{\bf SPM:} #1}}
\else  
  \newcommand{\pef}[1]{}
  \newcommand{\bda}[1]{}
  \newcommand{\rck}[1]{}
  \newcommand{\spm}[1]{}
\fi
\usetikzlibrary{calc,shapes,arrows,positioning,decorations.pathreplacing,
        intersections,matrix,patterns,fadings,calc,decorations.markings,arrows.meta}

\pgfplotsset{compat=1.18}
\headers{Automated Galerkin time stepping in \Irksome{}}{B.~D.~Andrews, P.~Brubeck, P.~E.~Farrell, R.~C.~Kirby, S.~P.~MacLachlan}
\title{Automated Galerkin time stepping in \Irksome{}\thanks{\funding{
BDA was supported by
the European Union [ERC, GeoFEM, 101164551].
PB and PEF were supported by
the Engineering and Physical Sciences Research Council [grant no.~EP/W026163/1],
the Science and Technology Facilities Council [grant no.~UKRI/ST/B000495/1],
and
the UKRI Digital Research Infrastructure Programme through the Science and Technology Facilities Council's Computational Science Centre for Research Communities (CoSeC).
PEF was supported by
the Donatio Universitatis Carolinae Chair ``Mathematical modelling of multicomponent systems'',
and
the National Science Foundation under grant no.~DMS-1929284 while in residence at the Institute for Computational and Experimental Research in Mathematics in Providence, USA.
The authors gratefully acknowledge the hospitality of the Erwin Schrödinger International Institute for Mathematics and Physics, Vienna, where part of this work was carried out.
The work of RCK was partially supported by National Science Foundation award 2410408.
The work of SPM was partially supported by an NSERC Discovery Grant.
For the purpose of open access, the authors have applied a CC BY public copyright license to any author accepted manuscript arising from this submission.
}}}
\author{Boris D.~Andrews\thanks{Mathematical Institute, University of Oxford, UK (\email{boris.andrews@maths.ox.ac.uk}).}
\and Pablo Brubeck\thanks{Mathematical Institute, University of Oxford, UK (\email{pablo.brubeck@maths.ox.ac.uk}).}
\and Patrick E.~Farrell\thanks{Mathematical Institute, University of Oxford, UK and Mathematical Institute, Faculty of Mathematics and Physics, Charles University, Czechia (\email{patrick.farrell@maths.ox.ac.uk}).}
\and Robert C.~Kirby\thanks{Department of Mathematics, Baylor University, USA (\email{robert\_kirby@baylor.edu}).}
\and Scott P.~MacLachlan\thanks{Department of Mathematics and Statistics, Memorial University of Newfoundland, Canada (\email{smaclachlan@mun.ca}).}}

\ifpdf
\hypersetup{
  pdftitle={Automated Galerkin time stepping in Irksome},
  pdfauthor={B. D. Andrews, P. Brubeck, P. E. Farrell, R. C. Kirby, and S. P. MacLachlan}
}
\fi

\begin{document}

\maketitle

\begin{abstract}
  As the study of temporal and spatial discretization schemes continues to advance, recent work has focused on the use of Galerkin-in-time discretization schemes that enable broader structure-preservation than is known for Runge--Kutta integrators.  While the promise of such discretizations is immense, their realization has, until now, generally relied on bespoke implementations that have limited their wider use.
In this work, we present automation in \Irksome{} for both discontinuous Galerkin and continuous Petrov--Galerkin time stepping of semidiscrete variational problems. The implementation supports auxiliary variables, flexible temporal quadrature, and monolithic algebraic solvers, and it enables switching between Runge--Kutta and Galerkin-in-time formulations with minimal changes to user code. Numerical examples illustrate accuracy, solver performance, and structure preservation across representative PDE systems.
\end{abstract}



\section{Introduction}\label{sec:introduction}

\pef{I'm not sure opening in this way is the right choice. Everyone knows to discretise ODEs you need to choose a method. Could we instead open by arguing that in most software this is hardcoded and laborious to change, whereas we're making it easy to change and experiment? In previous work we did this using the Runge--Kutta abstraction; this manuscript presents a major extension to other classes of time discretizations, with favourable properties?}
\rck{I agree with this approach -- we're pitching a paper on software rather than methods, although the methods are hard to come by in software.  I'm having a go at this here}
\bda{I've had a go too! :)}

Time stepping schemes are a critical aspect of the discretization of time-dependent partial differential equations.
Even when effective software like \Firedrake{}~\cite{FiredrakeUserManual} or FEniCS~\cite{BarattaEtal2023} streamlines the construction of spatial discretizations of time-dependent partial differential equations (PDEs), the choice of time discretization is typically hardcoded and laborious to change.
Further, methods of current interest, especially higher-order structure-preserving schemes, are not fully supported in general libraries for ordinary differential equations like Sundials~\cite{hindmarsh2005sundials} or PETSc/TS~\cite{abhyankar2018petsc}.
The \Irksome{}~\cite{farrell2021irksome, kirby2025extending, kirby2025automated} library took a first step to address this in \Firedrake{}:
by syntactic manipulation of the Unified Form Language (UFL)~\cite{Alnaes:2014}, the user need only supply a semidiscrete variational form, obtaining a broad range of Runge--Kutta (RK) and Runge--Kutta--Nystr\"om (RKN) schemes with essentially no further modifications.

Here, we extend this abstraction to Galerkin-in-time discretizations, supporting both discontinuous Galerkin (DG) and continuous Petrov--Galerkin (CPG) schemes.
The user-facing interface in \Irksome{} is unchanged;
switching between RK, DG, and CPG schemes typically requires changing only a single line of code.
Despite a well-developed theory, implementations of these schemes have until now generally been limited, cumbersome, or bespoke;
to our knowledge, this is the first general-purpose package to support them natively.
In contrast, while many schemes can be constructed in other frameworks, this often relies on ``repurposing'' spatial discretization tools that were not designed for time-integration, such as the implementation of DG methods in~\cite{SMAI-JCM_2023__9__61_0} within the DUNE framework~\cite{BASTIAN202175}.

A general semidiscrete form for PDE discretizations seeks
$u : \mathbb{R}_+ \to V_h$ such that
\begin{equation}\label{eq:problemg}
    G(t, u, \dot{u}; v) = 0,
\end{equation}
at all times $t$ and for all $v \in V_h$, where we write $\dot{u} := \tfrac{\partial}{\partial t}u$, $V_h$ is some (appropriately chosen) finite-element space, and $G$ is linear in $v$.
This formulation admits systems that are quasi-linear or potentially fully nonlinear in the time derivative;
its generality is central to the automation strategy in \Irksome{}, and the CPG methods we consider can be formulated and implemented at this level.
The DG schemes we consider, however, require more structure.
Thus, we also consider the sub-class of semidiscrete problems:
find $u : \mathbb{R}_+ \to V_h$ such that
\begin{equation}\label{eq:problemmf}
    M(\dot{u}, v) = F(t, u; v),
\end{equation}
at all times $t$ and for all $v \in V_h$, where $M$ is a bilinear form on $V_h$, while $F$ is linear in $v$.
For example, the heat equation with homogeneous Dirichlet or Neumann boundary conditions has a semidiscrete formulation
\begin{equation}\label{eq:heat}
    (\dot{u}, v) = (f, v) - (\nabla u, \nabla v),
\end{equation}
leading to $M(\dot{u}, v) = (\dot{u}, v)$ and $F(t, u; v) = (f, v) - (\nabla u, \nabla v)$.
The incompressible Navier--Stokes equations can also be written in the form of~\eqref{eq:problemmf} for vector-valued $u$ and $v$, where $M(\dot{u},v)$ would be a singular bilinear form omitting the time derivative of the pressure.

\pef{I wonder if we should first consider this general formulation, rather than introducing a more narrow one and immediately generalising?}
\rck{I agree and rearranged.}

\rck{Now, we need to edit/shorten the rest of the intro in light of this. We need to keep a bit of the history (e.g.~French) and recent work on auxiliary variables, etc, then pivot to the paper summary.}
\bda{I had a go at shortening!}

Although Runge--Kutta (RK) and related methods have long been applied to~\eqref{eq:problemmf}, they are not without limitations.
For example, the theory of geometric numerical integration~\cite{hairer2006geometric, marsden2001} shows they are unable, in general, to conserve non-quadratic invariants or dissipated quantities of interest~\cite{Celledoni_et_al_2009}.

\bda{Notes on why preserving such laws is important?}
\rck{Do we need this?  ``Everybody'' is interested in structure preservation.  My view is it's likely application-dependent and what your operational constraints are.  Preserving invariants is usually most important in long-time stability, but bounding the amount we can get them wrong by may be enough.  Other applications have hard constraints (e.g. positivity).  I don't think we need to say a lot here.}
\pef{I agree with Rob, I think most readers will appreciate why structure preservation matters.}

Finite elements in time, referred to interchangeably as variational-in-time and Galerkin-in-time discretizations, offer an alternative to \RK{} methods, whereby one interprets the dimension of time as one would interpret space in the construction of a finite-element discretization~\cite[Chap.~69~\&~70]{Ern_Guermond_2021c}.
Given the established connections between Galerkin-in-time and collocation schemes~\cite{huynh2011collocation}, one might ask why \Irksome{} (or any other time integration library) should separately support Galerkin-in-time methods when Runge--Kutta collocation support already exists~\cite{farrell2021irksome,kirby2025extending}.
The equivalence between Galerkin-in-time and Runge--Kutta methods relies on particular choices of temporal quadrature and, outside these choices, one obtains different schemes with potentially different structure preservation or other theoretical properties.
%% The current interest in Galerkin-in-time discretizations for structure preservation is of particular relevance in those cases where this equivalence does not hold.

\bda{More exposition on what FET is? Or is it assumed a reader of this manuscript would likely know that?}
\rck{I think we're fine if we cite the literature.  We're motivating them by saying they're like IRK but have different/sometimes better structure preservation and also describing the methods in the paper.}

It was first observed by French \& Schaeffer~\cite{French_Schaeffer_1990} that CPG discretizations could preserve non-quadratic conservation and dissipation laws for certain problems where RK methods do not.
These observations were later extended to general Hamiltonian systems written in canonical coordinates by Betsch \& Steinmann~\cite{Betsch_Steinmann_2000b,Betsch_Steinmann_2000a} and, in recent years, by Egger, Habrich, and Shashkov~\cite{Egger_Habrich_Shashkov_2021} to a broad class of conservative and dissipative PDEs.
Mixed continuous Petrov--Galerkin schemes with carefully chosen auxiliary variables have further been shown to be able to preserve multiple conservation and dissipation laws for many ODE~\cite{Andrews_Farrell_2025b,Cohen_Hairer_2011,Hairer_Lubich_2014,McLachlan_Quispel_Robidoux_1999} and PDE~\cite{andrews2025enforcing,Giesselmann_Karsai_Tscherpel_2025} systems.

Despite their advantages, however, general-purpose implementations of these schemes have remained unavailable.
They lead to very complicated algebraic systems that couple all of the temporal degrees of freedom within a time step, much like the stage coupling of fully implicit Runge--Kutta schemes.  Consequently, it is difficult to derive the system residual and Jacobian just from a black-box implementation of the ODE, and solving the equations requires special care. These implementation difficulties are magnified when auxiliary variables are introduced, as the auxiliary variables are sought in spaces with different continuity in time from the primal ones.

Here, we build on the \Irksome{} project~\cite{farrell2021irksome,kirby2025extending}, a time stepping library for Firedrake~\cite{FiredrakeUserManual} that automates general Runge--Kutta methods for spatially discretized partial differential equations.
We manipulate abstract syntax in UFL~\cite{Alnaes:2014} to derive the variational system to be solved at each time step.
Our implementation enables auxiliary variables for structure preservation in mixed and differential--algebraic systems and flexible quadrature in time, and we can leverage existing solver infrastructure~\cite{farrell2021pcpatch,kirby2018solver} to apply monolithic multigrid methods~\cite{abu2022monolithic,mmg} as well.

\pef{Is it too soon to mention FEniCSx? I think it broadens the appeal of the paper if it works beyond Firedrake.}
\rck{I mentioned fenics above in my redone first paragraph.  We should come back to it in conclusions/future work.  Pablo, Jorgen, and Ahsan are working on it}

To illustrate the interface, consider the Navier--Stokes equations posed on the unit square, with homogeneous initial conditions and lid-driven cavity boundary conditions on the velocity.
In the setup shown in \Cref{fig:codeldc}, we use a $16 \times 16$ triangular mesh, Taylor--Hood elements in space, and a continuous Petrov--Galerkin discretization in time for the velocity.
Here, we use a discontinuous discretization in time for the pressure (appropriate, as discussed below, because this is a differential-algebraic equation).  Simple keyword arguments are used to change this behavior, and quadrature and other options are controlled similarly.
Users specify the semidiscrete form in UFL, and \Irksome{} automatically constructs and solves the coupled variational problems at each time step.  The user code is a minimal modification to existing codes that solve a stationary problem or use \Irksome{} for Runge--Kutta methods with \Firedrake{}.

\begin{figure}
  \begin{lstlisting}[basicstyle={\ttfamily\footnotesize}]
from firedrake import *
from irksome import Dt, GalerkinCollocationScheme, TimeStepper

msh = UnitSquareMesh(16, 16)
V = VectorFunctionSpace(msh, 'CG', 2)
W = FunctionSpace(msh, 'CG', 1)
Z = V * W
up = Function(Z)
u, p = split(up)
v, q = TestFunctions(Z)

Re = Constant(10.0)
time_degree = 1

t = Constant(0.0)
dt = Constant(1.0 / 16)

F = (inner(Dt(u), v) * dx + inner(dot(grad(u), u), v) * dx
     + 1/Re * inner(grad(u), grad(v)) * dx - inner(p, div(v)) * dx
     - inner(div(u), q) * dx)

bcs = [DirichletBC(Z.sub(0), 0, (1, 2, 3)),
       DirichletBC(Z.sub(0), as_vector([1, 0]), (4,))]

scheme = GalerkinCollocationScheme(time_degree, quadrature_degree="auto")
stepper = TimeStepper(F, scheme, t, dt, up, bcs=bcs, aux_indices=[1])

while float(t) < 1.0:
    stepper.advance()
    t.assign(float(t) + float(dt))
  \end{lstlisting}
  \caption{\Firedrake/\Irksome~listing for the Navier--Stokes lid-driven cavity problem.}
  \label{fig:codeldc}
%  \Description{Code listing for the heat equation in Irksome.}
\end{figure}

The rest of the paper is organized as follows.
In \Cref{sec:method_formulation}, we present the formulation of DG and CPG time stepping methods, including treatment of auxiliary variables.
In \Cref{sec:automation}, we describe the UFL-based automation strategy implemented in \Irksome{}.
In \Cref{sec:algebraic_solvers}, we introduce the algebraic solver framework for the resulting coupled systems.
In \Cref{sec:numerical_results}, we report numerical results on various model problems.
In \Cref{sec:conclusions}, we conclude with discussion and future work.

\section{Method formulation}\label{sec:method_formulation}

\bda{I've gone through this section and tidied up some of the notation/explanations. Feel free to edit/remove my changes as necessary! :)}

In this section, we describe the DG and CPG time stepping schemes for semidiscrete variational problems~\eqref{eq:problemg} and~\eqref{eq:problemmf}, including the introduction of discontinuous-in-time auxiliary variables.
With our implementation in mind, we focus particularly on the expansion of these schemes in terms of bases for the finite-dimensional spaces, and on temporal quadrature.

\subsection{Notation}

We partition $[0, T]$ into $N_t$ time steps of size $\Delta t = \tfrac{T}{N_t}$ and set $t_n = n \Delta t$;
however, uniform time steps are purely a notational convenience, imposing no restriction on the time stepping methods we consider.
We refer to each interval by $I_n = (t_{n-1}, t_n)$ for $n \geq 1$, and assume that initial data, $u_{n-1} \approx u(t_{n-1})$, is known when we start computing on interval $I_n$.
When a function $u$ might have a discontinuity at $t_n$, we use the notation $u(t_{n,-})$ and $u(t_{n,+})$ to denote evaluating $u$ at $t_n$ from the left or right, respectively.

Over the interval $I_n$, denote the set of all polynomials of degree $s$ by $\mathcal{P}_s(I_n)$.
Note that this $s$ is only the degree of the temporal discretization and independent of whatever choices are made in the spatial discretization, leading to $V_h$.
To discretize the space of functions from $I_n$ to $V_h$, we denote by $V_{h,s}(I_n)$ the space of all polynomial maps from $I_n$ into $V_h$,
\begin{equation}\label{eq:vhk}
  V_{h,s}(I_n) := \mathcal{P}_s(I_n; V_h).
\end{equation}

\subsubsection*{Basis expansion \& quadrature}

Equipping $\mathcal{P}_s(I_n)$ with a basis $\{ \phi_j \}_{j=0}^s$, we can write elements of $V_{h,s}(I_n)$ as
\begin{subequations}\label{eq:basis}
\begin{equation}
    V_{h,s}(I_n)
    = \left\{\sum_{j=0}^s u_j \phi_j : u_j \in V_h\right\}.
\end{equation}
Further, for a basis $\{ \psi_i \}_{i=1}^{\dim V_h}$ of $V_h$, we can expand as
\begin{equation}
    V_{h,s}(I_n)
    = \left\{\sum_{i=1}^{\dim V_h} \sum_{j=0}^s u_{ij} \psi_i \phi_j : u_{ij} \in \mathbb{R}\right\}.
\end{equation}
Multi-index notation lets us express this more compactly.
Let $\mathbb{I}_{h,s} = \{1,\ldots,\dim V_h\} \times \{0,\ldots,s\}$. We can write members of $V_{h,s}(I_n)$ as
\begin{equation}
    V_{h,s}(I_n)
    = \left\{ \sum_{I \in \mathbb{I}_{h,s}} u_I \Psi_I : u_I \in \mathbb{R} \right\},
\end{equation}
\end{subequations}
where $\{\Psi_I := \psi_{I_1} \phi_{I_2}\}_{I \in \mathbb{I}_{h,s}}$ form a basis for $V_{h,s}(I_n)$.

\bda{There's a slight re-use of notation here, with the index $I \in \mathbb{I}_{h,s}$ and interval $I_n = (t_{n-1}, t_n)$. Probably fine.}

We introduce compact notation for the evaluation $u(t_*) \in V_h$ of $u \in V_{h,s}(I_n)$ at a particular time $t_* \in I_n$,
\begin{equation}
    u(t_*)
    = \sum_{j=0}^s u_j \phi_j(t_*)
    = \sum_{i=1}^{\dim V_h} \sum_{j=0}^s u_{ij} \psi_i \phi_j(t_*)
    = \sum_{I \in \mathbb{I}_{h,s}} u_I \Psi_I(t_*),
\end{equation}
where $\Psi_I(t_*) := \psi_{I_1} \phi_{I_2}(t_*) \in V_h$.
Practical implementation of this basis requires a total order on members $I = (I_1, I_2) \in \mathbb{I}_{h,s}$.
Two options present themselves:
(i)~One could index over $(\dim V_h) I_2 + I_1$ (letting the spatial index vary fastest) yielding a block $(s+1) \times (s+1)$ system, where each group of equations/unknowns are the coefficients of one of the $u_{I_2}$.
(ii)~One could alternatively index over $(s+1) I_1 + I_2$ (letting the time index vary fastest) yielding a block $\dim V_h \times \dim V_h$ system, collecting all coefficients that multiply a fixed spatial test function together.  Our implementation follows the first option, although we note that the second is more natural for some solvers.

\begin{remark}[Polynomial degree vs.~stage count]\label{rem:s_convention}
    We use $s$ to denote the polynomial degree in time of the solution, rather than the number of stages in the resultant discretization.
    For collocation \RK{} methods, these are equivalent;
    for DG, however, this is not the case.
    DG($s$) is an $(s+1)$-stage method, making DG(0) (i.e.~$s = 0$) the lowest-order DG method.
    The CPG scheme introduced below (\Cref{sec:cpg}) operates differently.
    Similar to typical collocation methods, strong enforcement of the initial data in CPG($s$) consumes one degree of freedom from the degree-$s$ polynomial space;
    consequently, the polynomial degree in time of the solution and the number of stages coincide at $s$, making CPG(1) (i.e.~$s = 1$) the lowest-order CPG method.
\end{remark}

Both the DG and CPG discretizations require integration over a space-time domain $\Omega \times I_n$ that, for general computation, requires appropriate quadrature.
Assuming such a scheme is already applied to the spatial component, we focus on the temporal integration, which we perform over a reference interval $(0, 1)$, using the standard pullback to generate $N_q$ quadrature points, $\xi_q \in I_n$, and weights, $w_q$ proportional to $\Delta t $, giving a rule of the form
\pef{If we're looking to cut space, I think the standard material below until the next subsection could be deleted.}
\bda{I'd agree!}
\rck{Commenting out.  We can put it back if needed}
%% \begin{subequations}
%% \begin{equation}
%%   \phi_i(t) = \hat{\phi}_i(\hat{t}),
%% \end{equation}
%% for $t = t^{n-1} + \hat{t} \, \Delta t$, with the chain rule giving
%% \begin{equation}
%%   \frac{\mathrm{d}\phi_i}{\mathrm{d}t}(t) = \frac{1}{\Delta t}\frac{\mathrm{d}\hat{\phi}_i}{\mathrm{d}\hat{t}}(\hat{t}).
%% \end{equation}
%% \end{subequations}
%% A reference quadrature rule over $(0, 1)$,
%% \begin{subequations}
%% \begin{equation}
%%   \int_0^1 \hat{f}(\hat{t})\,\mathrm{d}\hat{t} \approx \sum_{q=1}^{N_q} \hat{w}_q\, \hat{f}(\hat{\xi}_q),
%% \end{equation}
%% can be mapped to a quadrature rule on $I_n$ by
\begin{equation}\label{eq:time_quadrature}
  \int_{I_n} f(t)\,\mathrm{d}t \approx \sum_{q=1}^{N_q} w_q\, f(\xi_q),
\end{equation}
%% \end{subequations}
%% where $\xi_q = t^{n-1} + \hat{\xi}_q \, \Delta t$ and $w_q = \Delta t \, \hat{w}_q$.

\subsection{Discontinuous Galerkin (DG)}

In the DG approach, the left-hand side term $M(\dot{u}, v)$ in \eqref{eq:problemmf} is discretized in time over $I_n$ by a bilinear form $\mathcal{M}_n = \mathcal{M}_n^{(1)} + \mathcal{M}_n^{(2)}$ over $V_{h,s}(I_n)$, where
\begin{subequations}
\begin{align}
    \mathcal{M}_n^{(1)}(u, v)
    &:=
    \int_{I_n} M(\dot{u}(t), v(t)) \,\mathrm{d}t,  \label{eq:calm}\\
    \mathcal{M}_n^{(2)}(u, v)
    &:=
    M(u(t_{n-1,+}) - u_{n-1}, v(t_{n-1,+})),  \label{eq:calm2}
\end{align}
with $\mathcal{M}_n^{(2)}$ handling the discontinuity in time.
The right-hand side term $F(t, u; v)$ is similarly discretized through a (possibly nonlinear) operator $\mathcal{F}_n$ over $V_{h,s}(I_n)$
\begin{equation}\label{eq:calf}
    \mathcal{F}_n(u; v)
    := \int_{I_n} F(t, u(t); v(t)) \,\mathrm{d}t
\end{equation}
\end{subequations}
The DG discretization over $I_n$ is then:
Find $u \in V_{h,s}(I_n)$ such that
\begin{equation}
  \label{eq:dgvarprob}
  \mathcal{M}_n(u, v)
  \;(= \mathcal{M}_n^{(1)}(u, v) + \mathcal{M}_n^{(2)}(u, v))\;
  = \mathcal{F}_n(u; v)
\end{equation}
for all $v \in V_{h,s}(I_n)$.
After solving this variational problem, we find our subsequent data $u_n \in V_h$ at $t_n$ by evaluation
\begin{equation}
  \label{eq:dgupdate}
  u_n = u(t_{n,-}).
\end{equation}

Our formulation of DG in time is derived by integrating by parts on the time derivative term, replacing the value of $u$ at $t_n$ with the previous time value, and then integrating by parts again without the replacement.
This gives a ``strong'' form of the DG method.
An equivalent ``weak'' form of the method is obtained by omitting the second integration by parts.

\subsubsection*{Basis expansion \& quadrature}

The variational problem~\eqref{eq:dgvarprob} is naturally satisfied if and only if it holds for the basis $\{\Psi_I\}$ of $V_{h,s}(I_n)$.
Performing such an expansion and approximating each integral with an appropriate quadrature rule, we obtain algebraic equations for the coefficients of $\{u_I \in \mathbb{R}\}_{I \in \mathbb{I}_{h,s}}$ for $u = \sum_{I \in \mathbb{I}_{h,s}} u_I \Psi_I$.

%% To approximate the temporal integrals in~\eqref{eq:dgvarprob} with this quadrature rule, we introduce matrices:
%% \begin{equation}
%%   \label{eq:Vdef}
%%   \begin{split}
%%   V_{ij} & = \hat{\phi}_j(\hat{\xi}_i), \\
%%   V^w_{ij} & = \hat{w}_i V_{ij}
%%   \end{split}
%% \end{equation}
%% which tabulate the reference temporal basis functions at the quadrature points and scale by quadrature points, and
%% \begin{equation}
%%   \label{eq:Ddef}
%%     D_{ij} = \frac{\partial \hat{\phi}_j}{\partial \tau}(\hat{\xi}_i)
%% \end{equation}
%% which tabulates the derivatives of those basis functions.

%% Then, we also have that on each interval $I_n$,
%% \begin{equation}
%%   \phi_j(\xi_i) = \hat{\phi}_j(\hat{\xi}_i) = V_{ij}
%% \end{equation}
%% \begin{equation}
%%   \frac{\partial \phi_j}{\partial t}(\xi_i)
%%   = \frac{1}{\Delta t} \frac{\partial \hat{\phi}_j}{\partial \tau}(\hat{\xi}_i) = \frac{1}{\Delta t} D_{ij}
%% \end{equation}

%% With this notation, we can evaluate $U^n$ and its temporal derivative at each of the quadrature points:
%% \begin{equation}
%%   \label{eq:Untab}
%%   \begin{split}
%%   U^n(\xi_q) & = \sum_{\imath=0}^k U^n_\imath \phi_\imath(\xi_q)
%%   = \sum_{\imath=0}^k V_{q\imath} U^n_\imath \\
%%   \frac{\partial U^n}{\partial t}(\xi_q) & = \frac{1}{\Delta t}\sum_{\imath=0}^k D_{q\imath} U^n_{\imath}
%%   \end{split}
%% \end{equation}

Assuming that $M$, $F$ are approximated by some appropriate spatial quadrature rule, we may approximate the time integrals in $\mathcal{M}_n$, $\mathcal{F}_n$ using the quadrature rule \eqref{eq:time_quadrature}.
Considering first $\mathcal{M}_n^{(1)}$,
\begin{subequations}\label{eq:timeintM1quad}
\begin{align}
    \mathcal{M}_n^{(1)}(u, v)
        &\approx  \sum_{q=1}^{N_q} w_q M(\dot{u}(\xi_q), v(\xi_q))  \\
        &\approx  \sum_{q=1}^{N_q} w_q M\!\left(\sum_{J \in \mathbb{I}_{h,s}}u_J\dot{\Psi}_J(\xi_q), \sum_{I \in \mathbb{I}_{h,s}}v_I\Psi_I(\xi_q)\!\right)\!  \\
        &\approx  \!\sum_{I, J \in \mathbb{I}_{h,s}}\! v_I \!\left(\sum_{q=1}^{N_q} w_q M(\dot{\Psi}_J(\xi_q), \Psi_I(\xi_q))\!\right)\!u_J.
\end{align}
\end{subequations}
Let us denote the innermost sums in the linear term by $M_{I\!J}^{(1)}$, such that
\begin{equation}
    \mathcal{M}_n^{(1)}(u, v)  \approx  \!\sum_{I, J \in \mathbb{I}_{h,s}}\! v_I M_{I\!J}^{(1)} u_J.
\end{equation}
Using the expansion $\Psi_I := \psi_{I_1} \phi_{I_2}$, we can decompose $M_{I\!J}^{(1)}$ further as
%\begin{subequations}
\begin{equation}
    M_{I\!J}^{(1)}
        =  \sum_{q=1}^{N_q} w_q M(\psi_{J_1}\dot{\phi}_{J_2}(\xi_q), \psi_{I_1}\phi_{I_2}(\xi_q))
        =  M(\psi_{J_1}, \psi_{I_1}) \cdot \sum_{q=1}^{N_q} w_q \dot{\phi}_{J_2}(\xi_q)\phi_{I_2}(\xi_q).
\end{equation}
%\end{subequations}
Denote these terms by $M_{I_1\!J_1}^{(\psi)}$, $D_{I_2\!J_2}^{(\phi)}$ respectively.
We note that the matrix $M^{(1)}$, composed of $(M_{I\!J}^{(1)})$, can be written as the Kronecker product
\begin{equation}\label{eq:m1_kronecker}
    M^{(1)} =D^{(\phi)} \otimes  M^{(\psi)} ,
\end{equation}
where the $(s+1) \times (s+1)$ matrix $D^{(\phi)}$ and $\dim V_h \times \dim V_h$ matrix $M^{(\psi)}$  are composed of $(D_{ij}^{(\phi)})$ and $(M_{ij}^{(\psi)})$  respectively.
Note that, if the quadrature rule \eqref{eq:time_quadrature} is of a sufficiently high degree, we have exactly
\begin{equation}
  D_{ij}^{(\phi)} = \int_{I_n} \dot{\phi}_j \phi_i \mathrm{d}t.
\end{equation}

The second term $\mathcal{M}_n^{(2)}$ follows a similar decomposition,
\begin{equation}
    \mathcal{M}_n^{(2)}(u, v)
        \approx  \!\sum_{I, J \in \mathbb{I}_{h,s}}\! v_I M(\Psi_J(t_{n-1,+}), \Psi_I(t_{n-1,+})) u_J
        - \!\sum_{I \in \mathbb{I}_{h,s}} v_I M(u_{n-1}, \Psi_I(t_{n-1,+})).
\end{equation}
Denoting $M_{I\!J}^{(2)} := M(\Psi_J(t_{n-1,+}), \Psi_I(t_{n-1,+}))$, the linear term decomposes further as
\begin{equation}
    M_{I\!J}^{(2)}
        =  M(\psi_{J_1}, \psi_{I_1}) \cdot \phi_{J_2}(t_{n-1,+}) \phi_{I_2}(t_{n-1,+}).
\end{equation}
Denoting the latter term by $P_{I_2 J_2}^{(\phi)}$, the matrix $M^{(2)}$ composed of $(M_{I\!J}^{(2)})$ similarly satisfies a Kronecker decomposition,
\begin{equation}\label{eq:m2_kronecker}
    M^{(2)} =  P^{(\phi)} \otimes M^{(\psi)},
\end{equation}
with $M^{(\psi)}$ as above and $P^{(\phi)}$ the $(s+1) \times (s+1)$ matrix composed of $(P_{ij}^{(\phi)})$.
Note that, no matter the basis $\{ \phi_i \}$ for $\mathcal{P}_s(I_n)$, $P^{(\phi)}$ is a rank-1 matrix of the form $\mathbf{v} \mathbf{v}^\top$ with $\mathbf{v} = (\phi_{i}(t_{n-1,+}))$.

We lastly apply the quadrature rule to $\mathcal{F}_n$.
With the nonlinearity, little is gained by expanding $u$ as a sum;
we therefore restrict our attention to $v$:
\begin{equation}\label{eq:timeintFquad1}
    \mathcal{F}_n(u, v)
        \approx  \sum_{I\in \mathbb{I}_{h,s}} v_I \sum_{q=1}^{N_q} w_q F(\xi_q, u(\xi_q); \Psi_I(\xi_q)).
\end{equation}
Expanding $\Psi_I$ gives
\begin{equation}\label{eq:timeintFquad2}
    \sum_{q=1}^{N_q} w_q F(\xi_q, u(\xi_q); \Psi_I(\xi_q))
        =  \sum_{q=1}^{N_q} w_q F(\xi_q, u(\xi_q); \psi_{I_1}) \phi_{I_2}(\xi_q).
\end{equation}

After basis expansion and quadrature, the coefficients of $u$ are determined by the algebraic system
\begin{equation}
  \label{eq:dgalg}
\sum_{J \in \mathbb{I}_{h,s}} M_{IJ} u_{J} - \sum_{q=1}^{N_q} w_q F(\xi_q, u(\xi_q); \psi_{I_1}) \phi_{I_2}(\xi_q) = 0
\end{equation}
The equation corresponding to each test function requires evaluating $u$ at quadrature points, which couples together the coefficients corresponding to different temporal basis functions in much the same way as stages are coupled in fully implicit RK methods.

Linearization gives a block-structured Jacobian in terms of the Jacobian of $F$ evaluated at the temporal quadrature points.  
We also mention that appropriate choices of bases and quadrature rules can actually reproduce implicit Runge--Kutta schemes.
The connections are somewhat easier to make in the CPG context; here we point out that if one puts $\{ \xi_q \}_{q=0}^s$ as the Radau quadrature points and
$\{\phi_i\}_{i=0}^s$ as the Lagrange polynomials associated with those points, Equation~\eqref{eq:dgalg} greatly simplifies.
In this case, the sum over quadrature points collapses since $\phi_{I_2}(\xi_q) = \delta_{I_2, q}$, and one can show that the method exactly agrees with the
RadauIIA collocation method~\cite{huynh2011collocation}.

\subsubsection*{Splitting \& differing quadrature rules}

More generally, one might wish to apply different quadrature rules to different terms in the equation, especially those involving $F$, either to reduce costs or to ensure fidelity in conserved quantities.
Writing $F(t, u; v) = \sum_{i=1}^{n_F} F_i(t, u; v)$, the right-hand side term in \eqref{eq:dgvarprob} splits by linearity as $\mathcal{F}_n = \sum_{i=1}^{n_F} \mathcal{F}_{n,i}$, for
\begin{equation}
    \mathcal{F}_{n,i}(u; v)
        :=  \int_{I_n} F_i(t, u; v) \mathrm{d}t.
\end{equation}
Approximating each integral in this sum differently adds notational rather than conceptual difficulty, and is handled in our implementation.

\subsubsection*{Strongly-enforced boundary conditions}

Up to this point, the topic of strongly-enforced boundary conditions has remained unaddressed.
Suppose that~\eqref{eq:problemmf} carries an additional condition
\begin{equation}
    B u = g,
\end{equation}
where $B$ is some trace operator restricting members of $V_h$ to their values on a portion of the domain boundary, and $g = g(t)$ may be time-dependent.
We must then supplement~\eqref{eq:dgvarprob} with boundary conditions on $u \in V_{h,s}(I_n)$, or equivalently on each coefficient $u_j \in V_h$.

To this end, we form the $L^2(I_n; BV_h)$ projection of $g$ onto $\mathcal{P}_s(I_n; BV_h)$,
\begin{equation}
    \pi g = \sum_{i=0}^s g_i \phi_i,
\end{equation}
where each coefficient function $g_i \in BV_h$ lies in the same function space as $g$, the relevant trace of $V_h$.
The coefficients are determined by the linear system
\begin{equation}
    \sum_{j=0}^s M^{(\phi)}_{ij} g_j = \int_{I_n} g(t) \phi_i(t) \,\mathrm{d}t,
\end{equation}
where $M^{(\phi)}_{ij} := \int_{I_n} \phi_i \phi_j \,\mathrm{d}t$ and the right-hand side is approximated by numerical quadrature.
We then supplement~\eqref{eq:dgvarprob} with the boundary conditions
\begin{equation}\label{eq:strongBCs}
    B u_j = g_j, \qquad 0 \leq j \leq s.
\end{equation}

\subsection{Continuous Petrov--Galerkin (CPG)}\label{sec:cpg}

The CPG discretization over $I_n$ is largely similar to the DG discretization.
Two important distinctions from the DG method are that (i)~the test functions are sought in a lower-degree discontinuous space, and (ii)~the initial data is enforced strongly, implying the solution is continuous across time steps.
The absence of the jump term means that we may apply CPG schemes to equations of the form~\eqref{eq:problemg}.

The variational formulation of~\eqref{eq:problemg} is quite straightforward,
we seek $u \in V_{h,s}(I_n)$ such that $ u(t_{n-1})  = u_{n-1}$ and
\begin{equation}
  \label{eq:cpgvarG}
    \mathcal{G}_n(u, v) = 0 
\end{equation}
for all $v \in V_{h,s-1}(I_n)$, where
\begin{equation}
  \mathcal{G}_n(u, v) = \int_{I_n} G(t, u, \dot{u}; v) \,\mathrm{d}t.
\end{equation}

In the special case of~\eqref{eq:problemmf}, the continuity of $u$ and lack of a jump term means we simply define $\mathcal{M}_n := \mathcal{M}_n^{(1)}$
and seek $u \in V_{h,s}(I_n)$ with $ u(t_{n-1})  = u_{n-1}$ and
\begin{equation}\label{eq:cpgvarprob_pre}
  \mathcal{M}_n(u, v)
  = \mathcal{F}_n(u; v)
\end{equation}
for all $v \in V_{h,s-1}(I_n)$.

\subsubsection*{Lifting of the initial data}

To pose \eqref{eq:cpgvarprob_pre} over a non-affine space, define the restricted subspace $V_{h,s}^0 \subsetneq V_{h,s}$,
\begin{equation}\label{eq:vkh_restricted}
    V_{h,s}^0(I_n)
    := \{v \in V_{h,s} : v(t_{n-1}) = 0\}.
\end{equation}
We will impose the initial data in \eqref{eq:cpgvarprob_pre} through a lifting.

We continue to use $\{ \phi_j \}_{j=0}^{s-1}$ as the polynomial basis for the test space $\mathcal{P}_{s-1}(I_n)$ and let $\{ \varphi_i \}_{i=0}^{s}$ denote a basis for $\mathcal{P}_{s}(I_n)$ used in the trial space.
To work within the restricted subspace $V_{h,s}^0 \subsetneq V_{h,s}$ \eqref{eq:vkh_restricted}, we impose the condition on $\{ \varphi_j \}$ that
\begin{equation}\label{eq:left}
    \varphi_j(t_{n-1})
    = \begin{cases}
        1, & j = 0 \\
        0, & 1 \leq j \leq s
    \end{cases}.
\end{equation}
This property holds (among many other choices) for the Lagrange and Bernstein polynomials.  Then, $\{\varphi_i\}_{i=1}^s$ encodes a basis for the subspace vanishing at $t_{n-1}$.  

The CPG problem~\eqref{eq:cpgvarG} is then stated as
seeking $u \in u_{n-1}\varphi_0 + V_{h,s}^0(I_n)$ such that
\begin{equation}
  \mathcal{G}_n(u, v) = 0
\end{equation}
for all $v \in V_{h,s-1}(I_n)$ or, equivalently, for all $\{\phi_i\}_{i=0}^{s-1}$.  

Similarly, problem \eqref{eq:cpgvarprob_pre} may be recast as:
Find $u \in u_{n-1}\varphi_0 + V_{h,s}^0(I_n)$ such that
\begin{equation}\label{eq:cpgvarprob}
  \mathcal{M}_n(u, v)
  = \mathcal{F}_n(u; v)
\end{equation}
for all $v \in V_{h,s-1}(I_n)$.
In either case, we write $u = u_{n-1}\varphi_0 + \tilde{u}$ for $\tilde{u} \in V_{h,s}^0(I_n)$.

\subsubsection*{Basis expansion \& quadrature}

Multi-index notation again offers a compact way to work over a basis of $V_{h,s}^0$:
Letting $\mathbb{J}_{h,s} = \{1,\ldots,\dim V_h\} \times \{1,\ldots,s\}$, we can write members of $V_{h,s}^0(I_n)$ as
\begin{equation}
    V_{h,s}^0(I_n)
    = \left\{ \sum_{J \in \mathbb{J}_{h,s}} u_J \Psi_J : u_J \in \mathbb{R} \right\},
\end{equation}
where $\{\Psi_J := \psi_{J_1} \varphi_{J_2}\}_{J \in \mathbb{J}_{h,s}}$ form a basis for $V_{h,s}^0(I_n)$.

The basis expansion for $\mathcal{M}_n = \mathcal{M}_n^{(1)}$ is largely the same as in the DG case \eqref{eq:m1_kronecker},
\begin{equation}
    M  =  M^{(1)}  =  D^{(\phi)} \otimes M^{(\psi)},
\end{equation}
although the time component $D^{(\phi)} = (D_{ij}^{(\phi)})$ is redefined with entries
\begin{equation}
    D_{ij}^{(\phi)}  =  \sum_{q=1}^{N_q} w_q \dot{\varphi}_j(\xi_q) \phi_{i}(\xi_q).
\end{equation}

The connection between CPG and collocation-type schemes, as explored in~\cite{huynh2011collocation}, is clearer than for DG schemes.  Here, we choose an $s$-point quadrature rule such as Gauss--Legendre or right-Radau.  Then, we take the
test basis as the Lagrange polynomials of degree $s-1$ associated with these points.
Then, the integrated Lagrange polynomials form a basis for degree $s$ polynomials vanishing at the left endpoint.  Approximating all of the integrals with this scheme renders $D^{(\phi)}$ diagonal and exactly produces the collocation IRK scheme associated with these points.
This demonstrates both the overlap and difference between RK and CPG schemes. Not all implicit RK schemes arise from collocation, while Galerkin with higher-order quadrature schemes produces methods beyond Runge--Kutta.

%% \subsubsection*{Linear systems}

%% For linear systems, the basis expansion for $\mathcal{F}$ \eqref{eq:f_kronecker} is again largely similar,
%% \begin{equation}
%%     F  =  \sum_{q=1}^{N_q} w_q F_p^{(\psi)} \otimes K_p^{(\phi)},
%% \end{equation}
%% with the time component $K_p^{(\phi)} = (K_{ij,p}^{(\phi)})$ redefined for entries
%% \begin{equation}
%%     K_{ij,p}^{(\phi)}  =  \varphi_{j}(\xi_q)\phi_{i}(\xi_q).
%% \end{equation}
%% We obtain then a linear system very similar to~\eqref{eq:dglinsys_pre},
%% \begin{equation}
%%   \!\left[ M^{(\psi)} \otimes D^{(\phi)} - \sum_{q=1}^{N_q} w_q F_p^{(\psi)} \otimes K_p^{(\phi)} \right]\! \mathbf{u} = \mathbf{b}.
%% \end{equation}

%% For time-homogeneous linear systems with $\mathcal{F}(t; u, v)$ independent of $t$, the Kronecker structure for $F$ again simplifies as
%% \begin{equation}
%%     F  =  F^{(\psi)} \otimes K^{(\phi)},
%% \end{equation}
%% simplifying the linear system as
%% \begin{equation}
%%   \!\left[ M^{(\psi)} \otimes D^{(\phi)} - F^{(\psi)} \otimes K^{(\phi)} \right]\! \mathbf{u} = \mathbf{b}.
%% \end{equation}
%% Here, certain choices of bases and quadrature give exact equivalence to implicit RK schemes~\cite{huynh2011collocation}.
%% Also, like the DG formulation, right multiplication by $I_{\dim V_h} \otimes (D^{(\phi)} + P^{(\phi)})^{-1}$ implies a block structure with a traditional Butcher tableau and
%% left multiplication by $I_{\dim V_h} \otimes (K^{(\phi)})^{-1}$ implies a block structure mirroring that of Butcher~\cite{butcher1976implementation}.
%% \bda{Same comments here. :)}

\subsection{Auxiliary variables}\label{sec:aux_vars}

\pef{Boris and I should have a go at this at some point.}
\bda{I've re-written a little. :)}
\rck{We don't always need auxiliary variables to get these properties right.
  Sometimes Galerkin Just Works (see NLS and Allen--Cahn).  The intro text
  of this section should acknowledge that and say that we can often modify
  the system to DTRT when it doesn't.}
\bda{I've added a footnote. Is that enough do you think?}
\pef{I've rewritten it.}
\bda{I like that! :)}
Certain transient systems include variables on which a time derivative is not taken.
These include algebraic constraints in differential-algebraic equations, and additional variables in mixed formulations of PDEs.
For example, in a velocity-pressure formulation of the incompressible Navier--Stokes equations, the pressure has no time derivative; in the case
of a velocity-stress formulation, the stress has no time derivative.
We refer to these as auxiliary variables.
For many systems, CPG discretizations automatically preserve structural conservation laws or dissipation inequalities \cite{Egger_Habrich_Shashkov_2021}. If this is not the case, it is possible to introduce additional auxiliary variables to guarantee conservation or dissipation on discretization \cite{andrews2025enforcing}.

Whereas DG methods require no special treatment for such differential--algebraic systems,
a na\"ive CPG discretization would take these auxiliary variables to be in the degree-$s$ continuous-in-time solution space, imposing initial data strongly.
Since they are not differentiated in time, auxiliary variables require no initial conditions.
They may thus benefit from a different treatment when discretizing in time.

For some auxiliary space $W_h$, the general semidiscrete PDE formulation \eqref{eq:problemg} can be extended to accommodate such auxiliary variables:
Find $u: \mathbb{R}_+ \rightarrow V_h$ and $w: \mathbb{R}_+ \rightarrow W_h$ such that
\begin{equation}
    G_1(t, u, \dot{u}, w; v_1) = 0,  \qquad
    G_2(t, u, w; v_2) = 0,
\end{equation}
for all $v_1 \in V_h$ and $v_2 \in W_h$.
In the more restricted case of~\eqref{eq:problemmf}, we
find $u: \mathbb{R}^+ \rightarrow V_h$ and $w: \mathbb{R}^+ \rightarrow W_h$ such that
\begin{equation}\label{eq:problemmf_aux}
    M(\dot{u}, v_1) = F_1(t, u, w; v_1),  \qquad
    0 = F_2(t, u, w; v_2),
\end{equation}
at all times $t$ and for all $v_1 \in V_h$ and $v_2 \in W_h$.

% Here, $G_1$ rephrases the original PDE with the relevant expression replaced by $w$, and $G_2$ expresses the algebraic constraint between those expressions and the auxiliary variables.

In both the DG setting \eqref{eq:dgvarprob}, where solution variables are sought in a discontinuous-in-time space $V_{h,s}(I_n)$, and the CPG setting \eqref{eq:cpgvarprob}, where they are sought in a continuous-in-time space $V_{h,s}(I_n)$, we seek auxiliary variables $w_h$ in the discontinuous space
\begin{equation}\label{eq:whk}
  W_{h,s-1}(I_n) := \mathcal{P}_{s-1}(I_n; W_h).
\end{equation}
The auxiliary component of the system \eqref{eq:problemmf_aux} is then discretized over $I_n$ as
\begin{equation}
  0 = \mathcal{F}_{2,n}(u, w; v_2),
\end{equation}
for all $v_2 \in W_{h,s-1}$, where $\mathcal{F}_{2,n}$ is defined as in \eqref{eq:calf},
\begin{equation}\label{eq:calf2}
    \mathcal{F}_{2,n}(u, w; v)
    := \int_{I_n} F_2(t, u(t), w(t); v_2(t)) \,\mathrm{d}t.
\end{equation}
As with the original DG and CPG methods, we replace the time integration with summation using an arbitrary quadrature rule, and expand $w \in W_{h,s-1}(I_n)$ using a basis expansion as in the primal variables \eqref{eq:basis}.
 
By way of example, we study the CPG discretization of the incompressible Navier--Stokes equations, a differential-algebraic system.
Suppose $V_h$ and $W_h$ form a conforming inf-sup stable velocity--pressure pair.
Under appropriate boundary conditions, a suitable semidiscretization is:
Find $\mathbf{u} : \mathbb{R}_+ \to V_h$ and $p : \mathbb{R}_+ \to W_h$ such that
\begin{subequations}\label{eq:ns_semidiscrete}
\begin{equation}
    (\dot{\mathbf{u}}, \mathbf{v})  =  - \, (\mathbf{u} \cdot \nabla\mathbf{u}, \mathbf{v}) + (p, \mathrm{div}\,\mathbf{v}) - \frac{1}{\mathrm{Re}}(\nabla\mathbf{u}, \nabla\mathbf{v}),  \qquad
    0  =  (\mathrm{div}\,\mathbf{u}, q),
\end{equation}
at all times $t$ and for all $\mathbf{v} \in V_h$ and $q \in W_h$.
Writing the above system as
\begin{equation}
    M(\dot{\mathbf{u}}, \mathbf{v}) = F_1(\mathbf{u}, p; \mathbf{v}),  \qquad
    0  =  F_2(\mathbf{u}; q),
\end{equation}
\end{subequations}
we define $\mathcal{M}_n$, $\mathcal{F}_{1,n}$, $\mathcal{F}_{2,n}$ by integration of $M$, $F_1$, $F_2$ over $I_n$ as in (\ref{eq:calm},~\ref{eq:calf},~\ref{eq:calf2});
the space-time spaces $V_{h,s}$, $V_{h,s-1}$ and $W_{h,s-1}$ are similarly defined as in (\ref{eq:vhk},~\ref{eq:whk}).
A CPG discretization over $I_n$ of \eqref{eq:ns_semidiscrete} is then given as follows:
Find $\mathbf{u} \in V_{h,s}$ and $p \in W_{h,s-1}$ such that
\begin{subequations}
\begin{equation}
    \mathcal{M}_n(\mathbf{u}, \mathbf{v})  =  \mathcal{F}_{1,n}(\mathbf{u}, p; \mathbf{v}),  \qquad
    0  =  \mathcal{F}_{2,n}(\mathbf{u}; q),
\end{equation}
for all $\mathbf{v} \in V_{h,s-1}$ and $q \in W_{h,s-1}$, satisfying the initial condition
\begin{equation}
    \mathbf{u}(t_{n-1}) = \mathbf{u}_{n-1}.
\end{equation}
\end{subequations}
Note that $\mathbf{u}$ is sought in a continuous-in-time space with degree $s$ in time, while $p$, an auxiliary variable, is sought in a discontinuous-in-time space with degree $s-1$.

\section{Automation via UFL manipulation}\label{sec:automation}

In~\cite{farrell2021irksome}, an approach to UFL manipulation was introduced to obtain stage-coupled equations for \RK{} methods;
this has been extended to other \RK{} formulations in~\cite{kirby2025extending} and to \RKN{} methods in~\cite{kirby2025automated}.
Time discretizations are constructed through a UFL representation of the semidiscrete formulation, with time derivatives (e.g.~$\dot{u}$) represented through the \Irksome{} \texttt{Dt} function (e.g.~\texttt{Dt(u)}).
We keep the same interface here for describing variational problems in semidiscrete form, so that users can interchangeably apply either \RK{} or Galerkin-type methods.

While our exposition above sought to show the block structure, our implementation follows a different path.
The core principle is to absorb the quadrature weights onto the test function, via linearity.
This does not change the linear system, but it simplifies the UFL manipulation within \Irksome{}.
Fixing ideas by considering the DG case, instead of the manipulations in (\ref{eq:timeintM1quad}--\ref{eq:m1_kronecker}) for $\mathcal{M}_n^{(1)}$, we write
\begin{equation}
    \mathcal{M}_n^{(1)}(u, v)
        \;\left(:= \int_{I_n} M(\dot{u}(t), v(t)) \,\mathrm{d}t\right)\;
        \approx \sum_{q=1}^{N_q} M(\dot{u}(\xi_q), w_q v(\xi_q)).
\end{equation}
Expanding $w_q v(\xi_q) = w_q \sum_{i=0}^s v_i \phi_i(\xi_q)$ then gives the $s+1$ coupled systems as required.
A similar manipulation allows one to absorb the quadrature weights into the test function $v$ in the quadrature approximation of $\mathcal{F}_n$ (\ref{eq:timeintFquad1}--\ref{eq:timeintFquad2}).

% v_i \phi_j(\xi_q)

In a UFL representation of these equations, each ``stage'' $u_j$ of the unknown $u = \sum_{j=0}^s u_j \phi_j$ (in the DG case; $u = \sum_{j=0}^{s} u_j \varphi_j$ for CPG) will be represented by some member of a \lstinline{FunctionSpace} from \Firedrake, called \lstinline{V}.
For mixed systems, this space may be a product of an arbitrary number of \lstinline{FunctionSpace} instances, each of which may itself be a \lstinline{VectorFunctionSpace} or \lstinline{TensorFunctionSpace}.
The variational problem for the stages is posed on a new \lstinline{FunctionSpace} called \lstinline{Vbig}, a \lstinline{MixedFunctionSpace} comprising the $(s+1)$-way Cartesian product of the original \lstinline{V} with itself.
In addition to the UFL semidiscrete form, the \Irksome{} form manipulation for DG requires a univariate finite element basis for the temporal test and trial space and a quadrature rule, both of which are provided by FIAT~\cite{Kirby:2004}.

We then split up the UFL object for $u$ into a \lstinline{numpy} object array with $s+1$ rows, one for each $u_j$.
For mixed systems, when each $u_j$ comes from a mixed space, additional ranks of this array correspond to its components.
Left-multiplying this array by the matrix $(\phi_j(\xi_q))_{q,j}$ then gives a symbolic representation of the values $u(\xi_q)$.
We do the same pair of operations to the test function $v$, while also scaling each row by the quadrature weights $w_q$.
A similar operation tabulates the time derivative $\dot{u}$ at quadrature points $\xi_q$.

Quadrature approximation for $\mathcal{M}_n^{(1)}$ and $\mathcal{F}$ is achieved by looping over the (time) quadrature points, using UFL replacement to substitute (i)~$\xi_q$ for \lstinline{t}, (ii)~the time derivative at quadrature points for \lstinline{Dt(u)}, and (iii)~the test function values at quadrature points.
These substituted values are summed into a new UFL form.
We perform further UFL manipulation to evaluate $u(t_{n-1,+}) - u_{n-1}$ and $v(t_{n-1,+})$ for the term $\mathcal{M}_n^{(2)}$ \eqref{eq:calm2}.
A final loop sums over substitutions into the UFL expression for $\mathcal{F}_n$, and completes the construction of the discrete coupled variational problem.  This problem is augmented by any strong boundary conditions as specified in~\eqref{eq:strongBCs}.

\rck{TODO: Describe support for weak and strong forms of DG, it's just an option}

\section{Algebraic solvers}\label{sec:algebraic_solvers}

Sparse direct solvers are suitable for small model problems, including some considered in \Cref{sec:numerical_results} below.
However, larger problems require scalable solution algorithms.
Our generation of the UFL for the stage-coupled problems composes naturally with all of the extensive solver functionality offered by Firedrake \cite{kirby2018solver}.
While the algorithmic components are well-established, assembling them for the heavily coupled block structures of DG and CPG systems is typically a laborious, manual process. Our automation framework drastically simplifies this.
Here, as a demonstration, we adapt the stage-coupled monolithic multigrid schemes developed
in~\cite{abu2022monolithic,mmg,vanlent2005} for implicit Runge--Kutta methods to
DG time stepping schemes;
application to CPG is also possible with certain straightforward modifications.

For nonlinear problems, we first linearize by some Newton-type method, developing multigrid schemes for each linearized Newton iterate.
We thus consider general variational problems of the form
\begin{equation}\label{eq:general_linear}
  \mathcal{A}(u, v) = \mathcal{F}(v),
\end{equation}
where $\mathcal{A}$ is bilinear on $V_{h,s}(I_n) \times V_{h,s}(I_n)$.
Equivalently, expanding $u$ in the basis $\{\Psi_I\}$ and letting $v$ range over the same basis yields a general linear system
  $A\mathbf{u} = \mathbf{b}$.

The solver methods we consider are based on additive Schwarz relaxation schemes.
In the general steady-state setting, such methods decompose the spatial discretization space $V_h$ into a sum of possibly overlapping spaces
\begin{equation}\label{eq:schwarz}
  V_h = \sum_{i=1}^{N_S} V_i.
\end{equation}
When $V_h$ comprises continuous piecewise basis functions, a common decomposition is to define a subspace $V_i$ associated with each vertex $\mathbf{x}_i$ in the mesh, taking the span of all basis functions supported on the star of that vertex. %, as shown in \Cref{fig:starpatch}.
For certain symmetric and coercive problems, such a decomposition is known to give conditioning estimates independent of the polynomial degree~\cite{Pavarino:1993,Schoeberl:2008}, while independence is also observed empirically in far more general settings.
There is then a natural prolongation operator $p_i: V_i \rightarrow V_h$ defined simply by inclusion;
we encode its action by a matrix $P_i$.
We also provide a dual restriction operator $r_i: V_h^\prime \rightarrow V_i^\prime$;
this may be defined by transposition $R_i = (P_i)^\top$; however, other options are possible.
The additive Schwarz preconditioning operation $A \mapsto W^{-1}A$ is defined by specifying
\begin{equation}
  W^{-1} = \sum_{i=1}^{N_S} P_i A_i^{-1} R_i,
\end{equation}
where $A_i$ is the matrix corresponding to the linear system \eqref{eq:general_linear} restricted to $V_i$.
This amounts to solving $N_S$ decoupled, localized linear systems involving only the degrees of freedom in a patch of the spatial mesh.  

%% \begin{figure}[h]
%%   \centering
%%   \begin{tikzpicture}[scale=0.8]
%% \foreach \xa/\ya/\xb/\yb/\xc/\yc in {2/0/2/2/0/2, 4/0/4/2/2/2, 0/2/2/2/0/4, 2/2/4/2/2/4, 2/0/4/0/2/2, 2/2/2/4/0/4} {
%%     % Draw the triangle
%%     \draw[thin] (\xa,\ya) -- (\xb,\yb) -- (\xc,\yc) -- cycle;
%%     }

%% \foreach \i/\j in {2/1,3/1,1/2,2/2,3/2,1/3,2/3}{
%%   \draw[fill=black] (\i, \j) circle[radius=2pt];
%% }
    
%% \end{tikzpicture}
%% \caption{Typical vertex star patch for continuous piecewise quadratic basis functions.}
%% \label{fig:starpatch}
%% \end{figure}

In the transient setting, the general additive Schwarz decomposition \eqref{eq:schwarz} implies a decomposition of $V_{h,s}(I_n)$ by
\begin{equation}
V_{h,s}(I_n) = \sum_{i=1}^{N_S} V_i(I_n),
\end{equation}
where $V_i(I_n)$ comprises all polynomial maps from $I_n$ into $V_i$
\begin{equation}
V_i(I_n) = \left\{ \sum_{j=0}^s \phi_j u_j : u_j \in V_{i} \right\}.
\end{equation}
We may define corresponding prolongation and restriction operators for these spaces in terms of the given operations for $V_h$.
The prolongation $p_{i,s}: V_i(I_n) \rightarrow V_{h,s}(I_n)$ follows again by natural inclusion;
we can write it explicitly as
\begin{equation}
  p_{i,s}(u) = \sum_{j=0}^s \phi_j u_j = \sum_{j=0}^s \phi_j p_i(u_j).
\end{equation}
Ordering the basis coefficients $u_{ij}$ in a vector, with those of $u_j$ together, the operator $p_{i,s}$ can be represented as a matrix $P_{i,s} = I_{s+1} \otimes P_i$,
where $I_{s+1}$ is the $(s+1) \times (s+1)$ identity matrix.
With this ordering, we also have a restriction matrix
  $R_{i,s} = I_{s+1} \otimes R_i$;
in particular, when restriction is obtained by transposition, we have
  $R_{i,s} = P_{i,s}^\top = I_{s+1} \otimes P_i^\top$.
The preconditioning operation $A \mapsto W^{-1}A$ is then defined
\begin{equation}
  W^{-1} = \sum_{i=1}^{N_S} P_{i,s} A_{i,s}^{-1} R_{i,s},
\end{equation}
with $A_{i,s}$ defined similarly to $A_i$, by restriction of the linear system \eqref{eq:general_linear} to $V_{i,s}$.
The corresponding spatially localized linear systems are coupled across stages.

A geometric multigrid algorithm follows then in the usual way.
A sequence of nested spaces $V_{h_1} \subset V_{h_2} \subset \dots \subset V_{h_\ell} = V_{h}$  induces a sequence of nested spaces
$V_{h_1,s}(I_n) \subset \dots \subset V_{h_\ell,s}(I_n) = V_{h,s}(I_n)$.
Starting from the finest level, we perform some iteration (damped Richardson, Chebyshev, or a Krylov method) preconditioned with $W$, then transfer the residual from this iteration to the next coarser mesh.
On the coarsest mesh, we use a sparse direct solver or a Krylov method preconditioned with a processor-level additive Schwarz method using sparse LU or incomplete factorization.
Finally, we interpolate corrections from each mesh to the next finer mesh in the hierarchy, followed by additional preconditioned iteration on the resulting corrected approximation.
We report results for this solver strategy in \Cref{sec:incomp}.

\pef{Maybe give a solver diagram, if there's space?}

\section{Numerical results}\label{sec:numerical_results}

We demonstrate the new support within \Irksome{} on four representative PDE systems: CPG discretizations of Gross--Pitaevskii (a Hamiltonian system) and Allen--Cahn (a gradient-descent system), CPG and DG discretizations of incompressible flow, and a mixed CPG discretization of compressible flow.
As differential--algebraic schemes, the final two CPG discretizations necessitate the use of auxiliary variables (\Cref{sec:aux_vars}).
All source code is archived on Zenodo~\cite{zenodo/Zenodo-20260619.7}.

\bda{Add Zenodo citation/DoI}

\subsection{Hamiltonian system (Gross--Pitaevskii)}
CPG discretizations are known to conserve energy, independent of the time step size. For general Hamiltonian systems written in canonical coordinates~\cite{Betsch_Steinmann_2000b,Betsch_Steinmann_2000a},
this can readily be shown by taking the test function as the time derivative of the solution~\cite{Egger_Habrich_Shashkov_2021}.
As an illustrative example, we consider the defocusing nonlinear Schr\"odinger equation, otherwise known as the Gross--Pitaevskii equation~\cite{kevrekidis2015}, where the wave function satisfies
\begin{equation}
  i \psi_t = -\tfrac{1}{2} \Delta \psi + V(\mathbf{x}) \psi + \beta |\psi|^2 \psi,
\end{equation}
for $\mathbf{x} \in \Omega \subset \mathbb{R}^d$ with $d \in \{1, 2, 3\}$ and $\beta > 0$.
The function $V : \Omega \to \mathbb{R}$ represents a potential.
For simplicity, we consider Dirichlet boundary conditions $\psi = 0$ on $\partial \Omega$ and choose $\Omega$ large enough that this truncation has a negligible effect.

The equation has two important invariants:
(i) the wave function normalization,
\begin{subequations}
\begin{equation}
    N(\psi) = \int_\Omega |\psi(\mathbf{x}, t)|^2 \,\mathrm{d}x = 1, \ \ \ t \geq 0,
\end{equation}
and (ii) the energy,
\begin{equation}
    E(\psi) = \int_\Omega \frac{1}{2} \left| \nabla \psi(\mathbf{x}, t) \right|^2 + V(\mathbf{x}) \left|\psi(\mathbf{x}, t)\right|^2 + \frac{\beta}{2} \left| \psi(\mathbf{x}, t) \right|^4 \,\mathrm{d}x, \ \ \ t \geq 0.
\end{equation}
\end{subequations}

We compute the ground-state solution $\psi$ minimizing the energy $E(\psi)$ subject to $N(\psi) = 1$ via the techniques in~\cite{bao2004computing} and use this as the initial condition.
Many finite element and finite difference methods for the spatial discretization are proposed in the literature, and some time stepping schemes must pay special attention to the invariants.
We use a standard Galerkin spatial discretization with $P^1$ Lagrange elements. \Irksome{} allows us to conveniently compare Runge--Kutta and CPG time discretizations.

The implicit Gauss--Legendre (GL) methods, with the implicit midpoint rule being the lowest-order instance, are A-stable and symplectic, allowing the use of large time steps while preserving all quadratic invariants.
With this property, they exactly conserve $N$ (up to solver tolerances and roundoff) but do not exactly conserve the energy $E$.
In contrast, CPG exactly conserves $E$ (up to solver tolerances and roundoff) but does not maintain normalization $N$.
\Irksome{} allows switching between these methods with a single line of code.

\begin{figure}
  \begin{subfigure}[t]{0.475\textwidth}
\begin{tikzpicture}[scale=0.7]
\begin{axis}[
    xlabel = $t$,
    ylabel = {$N/N(0)$},
    scaled y ticks=false,
    yticklabel style={/pgf/number format/fixed, /pgf/number format/precision=4},
    legend pos=south west,
    legend style={font=\scriptsize, fill opacity=0.85, text opacity=1},
    % other axis options...
]

\pgfplotstableread[col sep=comma]{data/gp/gp_result_gl_1_T12.csv}\mytableglone
\pgfplotstablegetelem{0}{N}\of{\mytableglone}
\pgfmathsetmacro{\Nrefglone}{\pgfplotsretval}
\pgfplotstablecreatecol[create col/expr = {\thisrow{N} / \Nrefglone}]
                       {Nnormglone}\mytableglone
\pgfplotstableread[col sep=comma]{data/gp/gp_result_gl_2_T12.csv}\mytablegltwo
\pgfplotstablegetelem{0}{N}\of{\mytablegltwo}
\pgfmathsetmacro{\Nrefgltwo}{\pgfplotsretval}
\pgfplotstablecreatecol[create col/expr = {\thisrow{N} / \Nrefgltwo}]
                       {Nnormgltwo}\mytablegltwo
\pgfplotstableread[col sep=comma]{data/gp/gp_result_cpg_1_T12.csv}\mytablecpgone
\pgfplotstablegetelem{0}{N}\of{\mytablecpgone}
\pgfmathsetmacro{\Nrefcpgone}{\pgfplotsretval}
\pgfplotstablecreatecol[create col/expr = {\thisrow{N} / \Nrefcpgone}]
                       {Nnormcpgone}\mytablecpgone
\pgfplotstableread[col sep=comma]{data/gp/gp_result_cpg_2_T12.csv}\mytablecpgtwo
\pgfplotstablegetelem{0}{N}\of{\mytablecpgtwo}
\pgfmathsetmacro{\Nrefcpgtwo}{\pgfplotsretval}
\pgfplotstablecreatecol[create col/expr = {\thisrow{N} / \Nrefcpgtwo}]
                       {Nnormcpgtwo}\mytablecpgtwo                                              

\addplot[seabornblue, line width = 2, mark = x, mark size = 5, mark options = {solid}, mark repeat = 20, mark phase = 4]
table[x=t, y=Nnormglone] {\mytableglone};
\addlegendentry{GL(1)}
\addplot[seabornblue, line width = 2, dashed, mark = square, mark size = 3, mark options = {solid}, mark repeat = 20, mark phase = 8]
table[x=t, y=Nnormgltwo] {\mytablegltwo};
\addlegendentry{GL(2)}

\addplot[seabornred, line width = 2, mark = x, mark size = 5, mark options = {solid}, mark repeat = 20, mark phase = 12]
table[x=t, y=Nnormcpgone] {\mytablecpgone};
\addlegendentry{CPG(1)}
\addplot[seabornred, line width = 2, dashed, mark = square, mark size = 3, mark options = {solid}, mark repeat = 20, mark phase = 16]
table[x=t, y=Nnormcpgtwo] {\mytablecpgtwo};
\addlegendentry{CPG(2)}
\end{axis}

\end{tikzpicture}
  \end{subfigure}
  \begin{subfigure}[t]{0.475\textwidth}
    \begin{tikzpicture}[scale=0.7]
\begin{axis}[
    xlabel = $t$,
    ylabel = {$E/E(0)$},
    scaled y ticks=false,
    yticklabel style={/pgf/number format/fixed, /pgf/number format/precision=3},
    % other axis options...
  ]
 
\pgfplotstableread[col sep=comma]{data/gp/gp_result_gl_1_T12.csv}\mytableglone
\pgfplotstablegetelem{0}{E}\of{\mytableglone}
\pgfmathsetmacro{\Erefglone}{\pgfplotsretval}
\pgfplotstablecreatecol[create col/expr = {\thisrow{E} / \Erefglone}]
                       {Enormglone}\mytableglone
\pgfplotstableread[col sep=comma]{data/gp/gp_result_gl_2_T12.csv}\mytablegltwo
\pgfplotstablegetelem{0}{E}\of{\mytablegltwo}
\pgfmathsetmacro{\Erefgltwo}{\pgfplotsretval}
\pgfplotstablecreatecol[create col/expr = {\thisrow{E} / \Erefgltwo}]
                       {Enormgltwo}\mytablegltwo
\pgfplotstableread[col sep=comma]{data/gp/gp_result_cpg_1_T12.csv}\mytablecpgone
\pgfplotstablegetelem{0}{E}\of{\mytablecpgone}
\pgfmathsetmacro{\Erefcpgone}{\pgfplotsretval}
\pgfplotstablecreatecol[create col/expr = {\thisrow{E} / \Erefcpgone}]
                       {Enormcpgone}\mytablecpgone
\pgfplotstableread[col sep=comma]{data/gp/gp_result_cpg_2_T12.csv}\mytablecpgtwo
\pgfplotstablegetelem{0}{E}\of{\mytablecpgtwo}
\pgfmathsetmacro{\Erefcpgtwo}{\pgfplotsretval}
\pgfplotstablecreatecol[create col/expr = {\thisrow{E} / \Erefcpgtwo}]
                       {Enormcpgtwo}\mytablecpgtwo                                              

\addplot[seabornblue, line width = 2, mark = x, mark size = 5, mark options = {solid}, mark repeat = 20, mark phase = 4]
table[x=t, y=Enormglone] {\mytableglone};
\addplot[seabornblue, line width = 2, dashed, mark = square, mark size = 3, mark options = {solid}, mark repeat = 20, mark phase = 8]
table[x=t, y=Enormgltwo] {\mytablegltwo};

\addplot[seabornred, line width = 2, mark = x, mark size = 5, mark options = {solid}, mark repeat = 20, mark phase = 12]
table[x=t, y=Enormcpgone] {\mytablecpgone};
\addplot[seabornred, line width = 2, dashed, mark = square, mark size = 3, mark options = {solid}, mark repeat = 20, mark phase = 16]
table[x=t, y=Enormcpgtwo] {\mytablecpgtwo};
\end{axis}
\end{tikzpicture}  
\end{subfigure}
  \caption{Comparing conservation for the implicit Gauss--Legendre and CPG schemes for the Gross--Pitaevskii equation. The former preserves the normalization (shown at left), but not the energy (shown at right), while the reverse holds for CPG.}
  \label{fig:gp}
\end{figure}
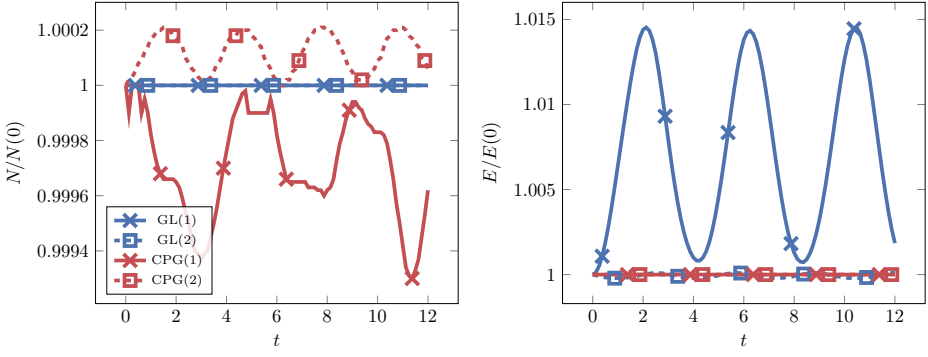

To illustrate the different behavior of the GL and CPG methods, we set $\Omega = (-8, 8)^2$, divided into a $32 \times 32$ mesh of squares, each subdivided into right triangles.
We use the potential
\begin{equation}
    V(x,y) = \frac{1}{2}\!\left(x^2 + y^2\right)\! + 4\exp(- (x - 1)^2 - y^2),
\end{equation}
and set $\beta = 200$.
This aligns with Case~II of the two-dimensional numerical experiments considered in the ground state calculations of~\cite{bao2004computing}.
We evolve the system using 1- and 2-stage GL and CPG schemes with $\Delta t = 0.125$, reporting the relative invariant conservation in \Cref{fig:gp}.
These computations confirm our expectations: the Gauss--Legendre schemes preserve $N$ but not $E$, while CPG preserves $E$ but not $N$.
%% Plots of the initial condition and the final state at time $t = 12$ for the 2-stage CPG simulation are shown in \Cref{fig:gpsol}.

\subsection{Gradient descent system (Allen--Cahn)}

A general gradient descent PDE in $u$ can be written in the variational form:
\begin{equation}
    (\dot{u}, v) + E'(u; v) = 0
\end{equation}
for all $v$, where $E'$ is the Fr\'echet derivative of some energy function $E(u) \in \mathbb{R}$ and $(\cdot, \cdot)$ is a general inner product. 
On the continuous level, taking $v = \dot{u}$ confirms that energy is dissipated:
\begin{equation}
    \frac{\mathrm{d}}{\mathrm{d}t}[E(u)]
        =  E'(u; \dot{u})
        =  - \, \| \dot{u} \|^2  \le  0.
\end{equation}
CPG methods are well-known to preserve this dissipation result:
\begin{equation}
    E|_{t_n} - E|_{t_{n-1}}
        = \int_{I_n} \frac{\mathrm{d}}{\mathrm{d}t}[E(u)]
        =  \int_{I_n} E'(u; \dot{u})
        =  - \, \int_{I_n} \| \dot{u} \|^2  \le  0,
\end{equation}
where we similarly consider $v = \dot{u}$~\cite{Egger_Habrich_Shashkov_2021}.
With exact integration, this holds unconditionally with respect to the time step size.

The situation is less clear for DG time stepping.
Selecting $v = \dot{u}$ leads to an energy relation with a sign-indefinite term.
Unconditional dissipativity in $L^2$ is shown in~\cite{estep2002dynamical} under some technical assumptions, but this does not prove energy monotonicity.
For a family of implicit Runge--Kutta methods (including Radau-IIA), Hairer and Lubich~\cite{Hairer_Lubich_2014} prove an energy-dissipation result under an abstract (if apparently mild) time step restriction.
Invariant energy quadratization~\cite{chen2025unconditionally} and scalar auxiliary variable methods~\cite{shen2018scalar} can also provide unconditional dissipation of a modified energy, typically in tandem with linearization, at the expense of additional unknowns.

The Allen--Cahn equation~\cite{allen1972,allen1979} models phase separation in multi-component alloys:
\begin{equation}\label{eq:acpde}
    \dot{u} - \epsilon^2 \Delta u + (u^3 - u)  =  0,
\end{equation}
posed in some domain $\Omega \subset \mathbb{R}^d$ for some diffusion parameter $\epsilon > 0$.
We impose homogeneous Neumann conditions on $\partial \Omega$.
Here the total free energy functional $E(u) \in \mathbb{R}$ is 
\begin{equation}
  E(u)  :=  \int_\Omega \frac{\epsilon^2}{2} |\nabla u|^2 + \frac{1}{2}(u^2-1)^2 \, \mathrm{d}x.
\end{equation}

We discretize~\eqref{eq:acpde} in space by a standard conforming Galerkin method with continuous piecewise linear basis functions on a $32 \times 32$ mesh of squares, subdivided into right triangles over the unit square domain.
The resulting semidiscrete variational problem is:
\begin{equation}
    (\dot{u}_h, v_h)
        + \epsilon^2(\nabla u_h, \nabla v_h)
        + (u_h^3 - u_h, v_h) = 0,
\end{equation}
where each inner product is over $L^2$.

\begin{figure}
  \begin{subfigure}[t]{0.475\textwidth}
    \caption{CPG}
    \begin{tikzpicture}[scale=0.7]
      \begin{axis}[xlabel={$t$}, ylabel={$E(t)$},
          ymode=log, xmode=log,
          legend style={font=\scriptsize, fill opacity=0.85, text opacity=1},
        legend pos = south west]
        \addplot[seabornblue, line width = 2, dashdotted, mark=x, mark size = 5, mark options = {solid}, mark repeat = 6] table[x=t,y=E, col sep=comma]{data/ac/energy.CPG.deg1.dt.0.1.csv};
        \addlegendentry{CPG(1) $\Delta t = 10^{-1}$}        
        \addplot[seaborngreen, line width = 2, dashed, mark=x, mark size = 5, mark options = {solid}, mark indices = {4, 9, 40, 90, 400}] table[x=t,y=E, col sep=comma]{data/ac/energy.CPG.deg1.dt.0.01.csv};
        \addlegendentry{CPG(1) $\Delta t = 10^{-2}$}
        \addplot[seabornred, line width = 2, dotted, mark=x, mark size = 5, mark options = {solid}, mark indices={2, 7, 20, 70, 200, 700, 2000}] table[x=t,y=E, col sep=comma]{data/ac/energy.CPG.deg1.dt.0.001.csv};
        \addlegendentry{CPG(1) $\Delta t = 10^{-3}$}        
        \addplot[seabornblue, line width = 2, dashdotted, mark=square, mark size = 3, mark options = {solid}, mark repeat = 6, mark phase = 2] table[x=t,y=E, col sep=comma]{data/ac/energy.CPG.deg2.dt.0.1.csv};
        \addlegendentry{CPG(2) $\Delta t = 10^{-1}$}        
        \addplot[seaborngreen, line width = 2, dashed, mark=square, mark size = 3, mark options = {solid}, mark indices = {4, 9, 40, 90, 400}] table[x=t,y=E, col sep=comma]{data/ac/energy.CPG.deg2.dt.0.01.csv};
        \addlegendentry{CPG(2) $\Delta t = 10^{-2}$}
        \addplot[seabornred, line width = 2, dotted, mark=square, mark size = 3, mark options = {solid}, mark indices={2, 7, 20, 70, 200, 700, 2000}] table[x=t,y=E, col sep=comma]{data/ac/energy.CPG.deg2.dt.0.001.csv};
        \addlegendentry{CPG(2) $\Delta t = 10^{-3}$}        
      \end{axis}
      \end{tikzpicture}
  \end{subfigure}
  \begin{subfigure}[t]{0.475\textwidth}
    \caption{DG}
    \begin{tikzpicture}[scale=0.7]
      \begin{axis}[xlabel={$t$}, ylabel={$E(t)$},
          xmode=log, ymode=log,
          legend style={font=\scriptsize, fill opacity=0.85, text opacity=1},
        legend pos = south west]
        \addplot[seabornblue, line width = 2, dashdotted, mark=x, mark size = 5, mark options = {solid}, mark repeat = 6] table[x=t,y=E, col sep=comma]{data/ac/energy.DG.deg1.dt.0.1.csv};
        \addlegendentry{DG(1) $\Delta t = 10^{-1}$}        
        \addplot[seaborngreen, line width = 2, dashed, mark=x, mark size = 5, mark options = {solid}, mark indices = {4, 9, 40, 90, 400}] table[x=t,y=E, col sep=comma]{data/ac/energy.DG.deg1.dt.0.01.csv};
        \addlegendentry{DG(1) $\Delta t = 10^{-2}$}
        \addplot[seabornred, line width = 2, dotted, mark=x, mark size = 5, mark options = {solid}, mark indices={2, 7, 20, 70, 200, 700, 2000}] table[x=t,y=E, col sep=comma]{data/ac/energy.DG.deg1.dt.0.001.csv};
        \addlegendentry{DG(1) $\Delta t = 10^{-3}$}        
        \addplot[seabornblue, line width = 2, dashdotted, mark=square, mark size = 3, mark options = {solid}, mark repeat = 6, mark phase = 2] table[x=t,y=E, col sep=comma]{data/ac/energy.DG.deg2.dt.0.1.csv};
        \addlegendentry{DG(2) $\Delta t = 10^{-1}$}        
        \addplot[seaborngreen, line width = 2, dashed, mark=square, mark size = 3, mark options = {solid}, mark indices = {4, 9, 40, 90, 400}] table[x=t,y=E, col sep=comma]{data/ac/energy.DG.deg2.dt.0.01.csv};
        \addlegendentry{DG(2) $\Delta t = 10^{-2}$}
        \addplot[seabornred, line width = 2, dotted, mark=square, mark size = 3, mark options = {solid}, mark indices={2, 7, 20, 70, 200, 700, 2000}] table[x=t,y=E, col sep=comma]{data/ac/energy.DG.deg2.dt.0.001.csv};
        \addlegendentry{DG(2) $\Delta t = 10^{-3}$}        
      \end{axis}
      \end{tikzpicture}
  \end{subfigure}  
  \caption{Monotonic energy decay for the Allen--Cahn equation discretized with CPG and DG schemes of degree 1 and 2 with varying time steps.}
  \label{fig:acenerg}
\end{figure}
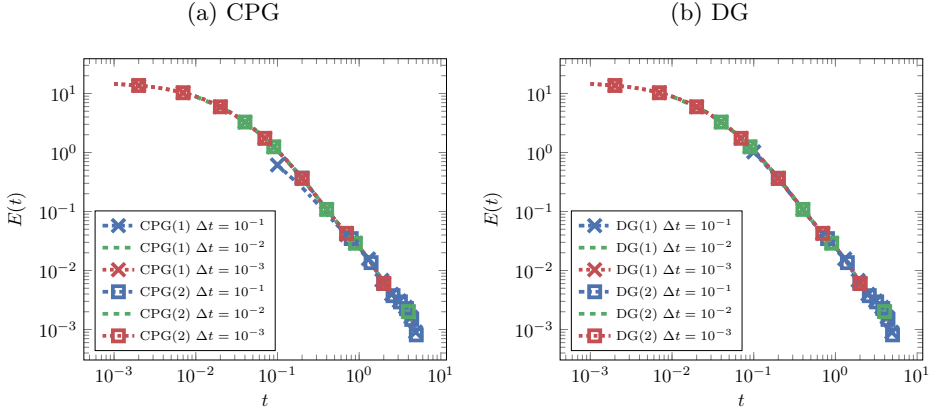

Here, we demonstrate energy dissipation using CPG(1), CPG(2), DG(1), and DG(2) methods for this problem.
We study a spinodal decomposition, taking a random initial condition and evolving in time.
We consider $\epsilon = 0.1$ and uniform time step sizes of $\Delta t \in \left\{10^{-3}, 10^{-2}, 10^{-1}\right\}$, integrating to final time $T=5.0$.
\Cref{fig:acenerg} shows monotonic energy decay over time.
We also observe monotonic energy decay for implicit Runge--Kutta schemes, but have not included these in the plots.
%% \Cref{fig:acsol} shows the initial condition and solution obtained using XXX.

\subsection{Incompressible flow (Navier--Stokes)}\label{sec:incomp}

Flow past a cylinder represents a frequently used benchmark problem for the 2D incompressible Navier--Stokes equations \cite{schafer1996benchmark,john2004reference}.
We consider both a CPG and DG discretization of this benchmark, the latter of which provides energy stability properties similar to those observed for the Radau-IIA methods in~\cite{abu2022monolithic}.
For time-homogeneous problems, the two methods are identical at lowest order $s=1$, where both reduce to implicit Euler (for fully linear systems, they are identical for all $s$).

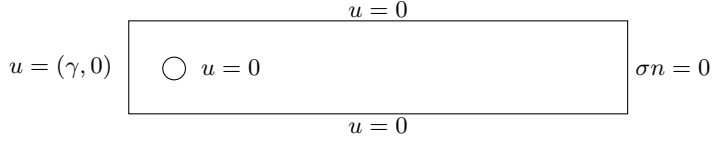
\begin{figure}
  \begin{center}
  \begin{tikzpicture}[scale=3]
    \draw (0,0) rectangle (2.2, 0.41);
    \draw (.2, .2) circle (0.05);
    \node[] at (1.1, 0.46) {\small $u=0$};
    \node[] at (1.1, -0.05) {\small $u=0$};
    \node[] at (.45, .205) {\small $u=0$};
    \node[] at (2.4, 0.205) {\small $\sigma n = 0$};
    \node[] at (-0.3, 0.205) {\small $u=\left(\gamma, 0\right)$};
  \end{tikzpicture}
  \end{center}
  \caption{Computational domain for Navier--Stokes flow past a cylinder.}
%  \Description{Picture showing computational domain for flow past a cylinder.}
  \label{fig:domain}
\end{figure}

The domain, shown in \Cref{fig:domain}, consists of the rectangle $(0, 2.2) \times (0, 0.41)$ with the circle of radius 0.05 centered at (0.2, 0.2) omitted.
On the left edge, we impose an inflow condition, setting the horizontal velocity component to be
\begin{equation}\label{eq:incomp_bc}
    \gamma(y, t) = 6 \sin \!\left(\frac{\pi t}{8} \right)\! \cdot \frac{y(0.41-y)}{0.41^2}
\end{equation}
and the vertical component to be zero.
Along the top and bottom edges and boundary of the circle we impose no-slip conditions.
We use a viscosity $\nu = 10^{-3}$, and integrate to final time $T = 8$, measuring the drag and lift on the circle boundary at each time step.
This aligns with the setting considered in \cite{john2004reference}, with a minor difference in how the outlet is handled:
while we consider natural traction-free outflow conditions, John imposes identical Dirichlet boundary conditions on the right edge as the left \eqref{eq:incomp_bc}.
Given the distance of the outlet from the obstacle, we expect this distinction to have a minimal effect on the computed drag and lift forces.

Two spatial discretizations (on quadrilateral and triangular grids) are considered in \cite{john2004reference}, using both a Crank--Nicolson temporal discretization and a certain fractional-step $\theta$-scheme.
The finest computation reported used a quadrilateral mesh with approximately 53k elements;
velocities were discretized with continuous biquadratic elements and pressures with discontinuous linear elements.
Results are presented for time step size from $0.04$ to $0.00125$.

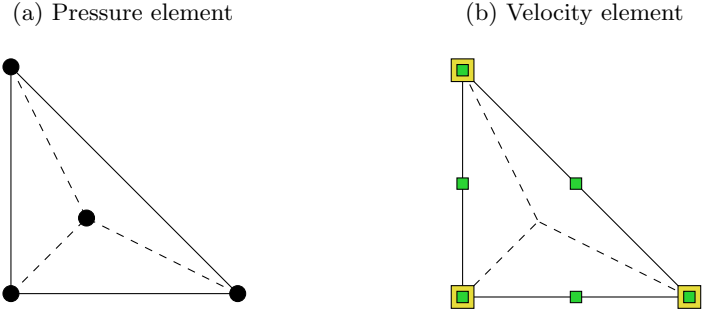
\begin{figure}
  \centering
  \begin{subfigure}[t]{0.45\textwidth}
    \centering
    \caption{Pressure element}
    \label{fig:aspressure}
    \begin{tikzpicture}[scale=1.5]
      \draw[thin] (0,0) -- (2,0) -- (0,2) -- cycle;
      \draw[thin, dashed] (0,0) -- (2/3,2/3);
      \draw[thin, dashed] (2,0) -- (2/3,2/3);
      \draw[thin, dashed] (0,2) -- (2/3,2/3);

      % Pressure dofs
      \draw[fill=black] (0,0) circle[radius=2pt];
      \draw[fill=black] (2,0) circle[radius=2pt];
      \draw[fill=black] (0,2) circle[radius=2pt];
      \draw[fill=black] (2/3,2/3) circle[radius=2pt];
    \end{tikzpicture}
  \end{subfigure}
  \begin{subfigure}[t]{0.45\textwidth}
    \centering
    \caption{Velocity element}
    \label{fig:asvelocity}
    \begin{tikzpicture}[scale=1.5]
      \draw[thin] (0,0) -- (2,0) -- (0,2) -- cycle;
      \draw[thin, dashed] (0,0) -- (2/3,2/3);
      \draw[thin, dashed] (2,0) -- (2/3,2/3);
      \draw[thin, dashed] (0,2) -- (2/3,2/3);
      \foreach \i/\j in {0/0, 2/0, 0/2} {
        \draw[fill=seabornsand!50!yellow] (\i-0.1, \j-0.1) rectangle (\i+0.1, \j+0.1);
      }

      \foreach \i/\j in {0/0, 1/0, 2/0, 0/1, 1/1, 0/2} {
        \draw[fill=seaborngreen!50!green] (\i-0.05, \j-0.05) rectangle (\i+0.05, \j+0.05);
      }
    \end{tikzpicture}
  \end{subfigure}
  \caption{Schematic showing the degrees of freedom for pressure (left) and velocity (right) for the Alfeld--Sorokina Stokes pair. Pressure degrees of freedom are shown with black dots.  Velocity values (as a pair of vector components at a point) are shown with green squares, while pointwise divergence of velocity is shown with yellow squares.}
  \label{fig:asdofs}
\end{figure}

We apply a different discretization to this benchmark.
Using Firedrake's interface to Netgen~\cite{betteridge2024ngspetsc}, we generate a hierarchy of triangular meshes whose edges interpolate the circle.  
The finest such mesh contains approximately 29k vertices and 57k triangles.
For the spatial discretization, we use the degree-2 Alfeld--Sorokina element for the velocity~\cite{alfeld-sorokina16}, available in Firedrake through our work on macroelements~\cite{brubeck2025fiat}.
This is a macroelement, defined over a barycentrically refined mesh;
taking the pressure in the degree-1 Lagrange space over the refined mesh (realized also as a macroelement in \Firedrake) and the velocity as a piecewise quadratic function over the original mesh with continuous divergence gives an inf-sup stable pair with a pointwise divergence-free velocity, whose degrees of freedom are illustrated in \Cref{fig:asdofs}.
The divergence-free kernel within the Alfeld--Sorokina space aligns exactly with that of the Scott--Vogelius pair of the same degree over the same split mesh~\cite{guzman2018inf, scott1984conforming}; however, the Alfeld--Sorokina velocity and pressure are supersmooth ($H^1(\mathrm{div})$- and $H^1$-conforming respectively) leading to much smaller spaces.
To maintain the exact pointwise divergence-free property without requiring complex Piola transformations, the curved boundary of the cylinder is approximated here by a polygonal mesh comprised of straight-edged affine triangles.

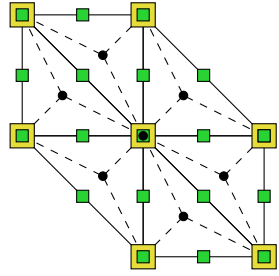
\begin{figure}[h]
  \centering
  \begin{tikzpicture}[scale=0.8]
\foreach \xa/\ya/\xb/\yb/\xc/\yc in {2/0/2/2/0/2, 4/0/4/2/2/2, 0/2/2/2/0/4, 2/2/4/2/2/4, 2/0/4/0/2/2, 2/2/2/4/0/4} {
    % Draw the triangle
    \draw[thin] (\xa,\ya) -- (\xb,\yb) -- (\xc,\yc) -- cycle;
    
    % Compute centroid
    \pgfmathsetmacro{\centx}{(\xa + \xb + \xc)/3}
    \pgfmathsetmacro{\centy}{(\ya + \yb + \yc)/3}
    
    % Draw dashed lines from vertices to centroid
    \draw[thin, dashed] (\xa,\ya) -- (\centx, \centy);
    \draw[thin, dashed] (\xb,\yb) -- (\centx, \centy);
    \draw[thin, dashed] (\xc,\yc) -- (\centx, \centy);

    \draw[fill=black] (\centx, \centy) circle[radius=2pt];
}

\foreach \i/\j in {2/0,4/0,0/2,2/2,4/2,0/4,2/4} {
  \draw[fill=seabornsand!50!yellow] (\i-0.2, \j-0.2) rectangle (\i+0.2, \j+0.2);
}

\foreach \i/\j in {2/0,4/0,0/2,2/2,4/2,0/4,2/4,3/0,1/1,2/1,3/1,4/1,1/2,3/2,4/2,0/3,1/3,2/3,3/3,1/4} {
  \draw[fill=seaborngreen!50!green] (\i-0.1, \j-0.1) rectangle (\i+0.1, \j+0.1);
}
\draw[fill=black] (2, 2) circle[radius=2pt];
    
\end{tikzpicture}
\caption{Typical Vanka patch for Alfeld--Sorokina pair.  Pressure degrees of freedom (omitted on the boundary) are shown with black dots.  Velocity values (as a pair of vector components at a point) are shown with green squares, while pointwise divergence of velocity is shown with yellow squares.}
\label{fig:vankanse}
\end{figure}

At each time step, we solve the stage-coupled variational problem using a Newton--Krylov method.
For each Newton step, we use GMRES preconditioned by a monolithic multigrid scheme as described in \Cref{sec:algebraic_solvers}.
The relaxation is based on a Vanka-type decomposition~\cite{vanka1986block,abu2022monolithic} similar to that used for higher-order Taylor--Hood elements in~\cite{rafiei2025improvingpatchselectionmonolithic}.
We show the single-stage decomposition in \Cref{fig:vankanse}.

\begin{figure}
  \begin{subfigure}[t]{0.45\textwidth}
    \caption{Drag}
  \begin{tikzpicture}[scale=0.7]
    \begin{axis}[xlabel={$t$}, ylabel={$C_D$}]
      \addplot[mark=none, line width = 1, each nth point=16, filter discard warning=false, restrict x to domain=0.001:10]
      table [x=t,y=CD,col sep=comma]{data/cylinder/dl_lvl4_dg2_Nt3200.dat};
    \end{axis}
  \end{tikzpicture}
  \end{subfigure}
  \begin{subfigure}[t]{0.45\textwidth}
    \caption{Lift}
  \begin{tikzpicture}[scale=0.7]
    \begin{axis}[xlabel={$t$}, ylabel={$C_L$}]
      \addplot[mark=none, line width = 1, each nth point=2, filter discard warning=false, restrict x to domain=0.001:10]
      table [x=t,y=CL,col sep=comma]{data/cylinder/dl_lvl4_dg2_Nt3200.dat};
    \end{axis}
  \end{tikzpicture}
  \end{subfigure}
  \caption{Drag and lift versus time for Navier--Stokes flow past a cylinder using $\Delta t = 0.0025$ and a DG(2) discretization in time.}
  \label{fig:dl}
\end{figure}
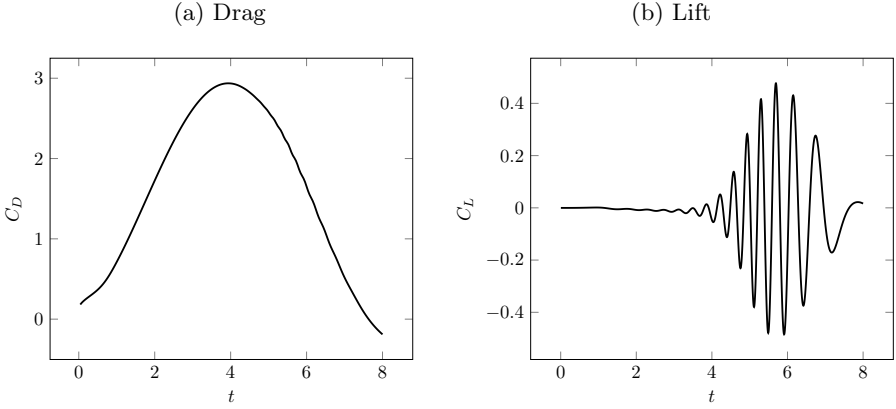

We first establish that our highest-fidelity DG computation (DG(2) with 3200 time steps, i.e.~$\Delta t = 0.0025$) can reproduce the benchmark lift and drag data reported in \cite{john2004reference}.
The differences in the discretization and outlet condition mean that exact agreement is not to be expected;
this initial comparison is intended only to confirm broad consistency with the expected dynamics.
% The convergence study that follows takes our own high-resolution computation as the reference solution.
Denoting by $\epsilon_D$ and $\epsilon_L$ the differences in drag and lift between our computation and the tabulated values, we measure their magnitudes in both $L^2([0,8])$ and at the final time $t = 8$.
To approximate the time integrals, we fit cubic splines to the tabulated data in \cite[Fig.~4]{john2004reference}, and integrate using \lstinline{scipy}.
\Cref{table:cfjohn} shows these errors are relatively small, consistent with the methodological differences noted above.
A plot of the lift and drag from our solution is shown in \Cref{fig:dl}.

\begin{table}
  \begin{center}
  \begin{tabular}{l|c}
    Quantity & Relative error \\ \hline
    $\| \epsilon_D \|_{L^2}$ & $2.53 \times 10^{-3}$ \\
    % $\| \epsilon_L \|_{L^2}$ & $6.85 \times 10^{-3}$ \\
    $|\epsilon_D(8)|$ & $2.51 \times 10^{-3}$ \\
    % $|\epsilon_L(8)|$ & $3.12 \times 10^{-2}$
  \end{tabular}
  \begin{tabular}{l|c}
    Quantity & Relative error \\ \hline
    % $\| \epsilon_D \|_{L^2}$ & $2.53 \times 10^{-3}$ \\
    $\| \epsilon_L \|_{L^2}$ & $6.85 \times 10^{-3}$ \\
    % $|\epsilon_D(8)|$ & $2.51 \times 10^{-3}$ \\
    $|\epsilon_L(8)|$ & $3.12 \times 10^{-2}$
  \end{tabular}
  \end{center}
  \caption{Comparison of lift and drag for Navier--Stokes flow past a cylinder using the Alfeld--Sorokina Stokes pair on a mesh with 29k vertices and DG(2) temporal discretization with $\Delta t = 0.0025$ against the finest computation reported in~\cite{john2004reference}.}
  \label{table:cfjohn}
\end{table}

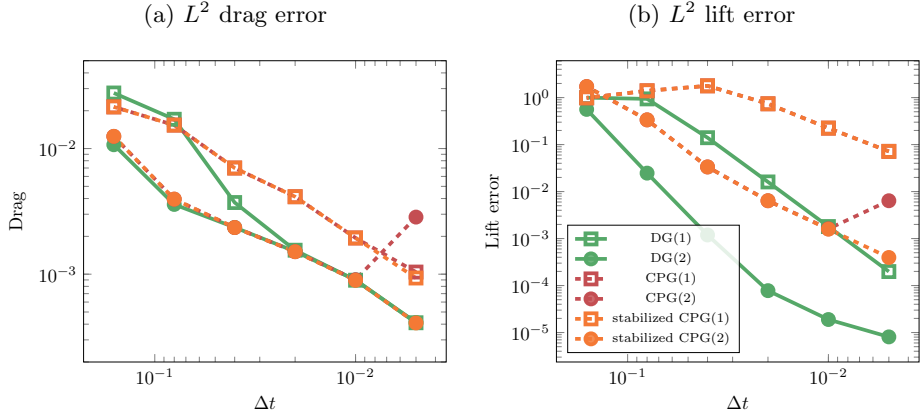
\begin{figure}
  \begin{subfigure}[t]{0.475\textwidth}
    \caption{$L^2$ drag error}
    \begin{tikzpicture}[scale=0.7]
      \begin{loglogaxis}[
          xlabel={$\Delta t$}, ylabel={Drag},
          x dir = reverse,
          ylabel near ticks,
          ymax=5.0e-2, ymin=2.e-4,
          % legend pos=south west,
          % legend style={font=\scriptsize, fill opacity=0.85, text opacity=1},
          % cycle list name=color list
          ]
        \addplot[seaborngreen, line width = 2, mark=square, mark size = 3]
        table [x=dt,y=L2drag, col sep=comma]{data/cylinder/dl_err_dt_dg1.csv};
        % \addlegendentry{DG(1)}

        \addplot[seaborngreen, line width = 2, mark=*, mark size = 3]
        table [x=dt,y=L2drag, col sep=comma]{data/cylinder/dl_err_dt_dg2.csv};
        % \addlegendentry{DG(2)}

        \addplot[seabornred, line width = 2, mark=square, mark size = 3, mark options = {solid}, dashed]
        table [x=dt,y=L2drag, col sep=comma]{data/cylinder/dl_err_dt_cpgaux1.csv};
        % \addlegendentry{CPG(1)}

        \addplot[seabornred, line width = 2, mark=*, mark size = 3, mark options = {solid}, dashed]
        table [x=dt,y=L2drag, col sep=comma]{data/cylinder/dl_err_dt_cpgaux2.csv};
        % \addlegendentry{CPG(2)}
        
        \addplot[seabornorange, line width = 2, mark=square, mark size = 3, mark options = {solid}, dashed]
        table [x=dt,y=L2drag, col sep=comma]{data/cylinder/dl_err_dt_cpg_aux_stable1.csv};
        % \addlegendentry{stabilized CPG(1)}

        \addplot[seabornorange, line width = 2, mark=*, mark size = 3, mark options = {solid}, dashed]
        table [x=dt,y=L2drag, col sep=comma]{data/cylinder/dl_err_dt_cpg_aux_stable2.csv};
        % \addlegendentry{stabilized CPG(2)}
        
      \end{loglogaxis}
    \end{tikzpicture}
  \end{subfigure}
  \begin{subfigure}[t]{0.475\textwidth}
    \caption{$L^2$ lift error}
    \begin{tikzpicture}[scale=0.7]
      \begin{loglogaxis}[
          xlabel={$\Delta t$}, ylabel={Lift error},
          x dir=reverse,
          ylabel near ticks,
          % ymax=1.e0, ymin=1.e-4,
          legend pos=south west,
          legend style={font=\scriptsize, fill opacity=0.85, text opacity=1},
          cycle list name=color list
          ]
        \addplot[seaborngreen, line width = 2, mark=square, mark size = 3]
        table [x=dt,y=L2lift, col sep=comma]{data/cylinder/dl_err_dt_dg1.csv};
        \addlegendentry{DG(1)}

        \addplot[seaborngreen, line width = 2, mark=*, mark size = 3]
        table [x=dt,y=L2lift, col sep=comma]{data/cylinder/dl_err_dt_dg2.csv};
        \addlegendentry{DG(2)}

        \addplot[seabornred, line width = 2, mark=square, mark size = 3, mark options = {solid}, dashed]
        table [x=dt,y=L2lift, col sep=comma]{data/cylinder/dl_err_dt_cpgaux1.csv};
        \addlegendentry{CPG(1)}
        
        \addplot[seabornred, line width = 2, mark=*, mark size = 3, mark options = {solid}, dashed]
        table [x=dt,y=L2lift, col sep=comma]{data/cylinder/dl_err_dt_cpgaux2.csv};
        \addlegendentry{CPG(2)}

        \addplot[seabornorange, line width = 2, mark=square, mark size = 3, mark options = {solid}, dashed]
        table [x=dt,y=L2lift, col sep=comma]{data/cylinder/dl_err_dt_cpg_aux_stable1.csv};
        \addlegendentry{stabilized CPG(1)}
        
        \addplot[seabornorange, line width = 2, mark=*, mark size = 3, mark options = {solid}, dashed]
        table [x=dt,y=L2lift, col sep=comma]{data/cylinder/dl_err_dt_cpg_aux_stable2.csv};
        \addlegendentry{stabilized CPG(2)}        
      \end{loglogaxis}
    \end{tikzpicture}
  \end{subfigure}
  \caption{Relative $L^2([0, 8])$ error in lift and drag versus the time step size for DG and CPG time stepping schemes. ``CPG'' refers to CPG with the pressure treated as a discontinuous-in-time auxiliary variable; ``stabilized CPG'' additionally under-integrates the $(\mathrm{div}\,\mathbf{u}_h, q_h)$ term using $s$-stage right-Radau quadrature.}
  \label{fig:dlerrdt}
\end{figure}

We now apply both DG and CPG time discretizations to this spatial discretization, gradually decreasing $\Delta t$ from $0.16$ to $0.005$ and comparing the lift and drag against a high resolution DG(2) reference solution with $\Delta t = 0.0025$.
These results appear in \Cref{fig:dlerrdt}.
For instance, to maintain both quantities from the DG solution within 1\% of the reference solution, we require $\Delta t \le 0.01$ for DG(1) and $\Delta t \le 0.04$ for DG(2).
The CPG discretization requires more care.
With no special treatment of the pressure, the discretization fails to converge to the reference solution;
for CPG(2) for instance, the $L^2$ error in the drag remains uniformly above $0.02$ for all $\Delta t$.
Convergence is restored by treating the pressure as a discontinuous-in-time auxiliary variable (see \Cref{sec:aux_vars}).
\Cref{fig:dlerrdt} shows the resulting CPG scheme converging at the expected rates for moderate $\Delta t$.

Despite this, CPG(2) shows an uptick in both drag and lift error at the finest time step $\Delta t = 0.005$.
This can be attributed to the behavior of the incompressibility constraint $\mathrm{div}\,\mathbf{u}_h = 0$.
With a standard CPG discretization, the error in this constraint at a given time $t_{n+1}$ can be bounded by the error at $t_n$ plus solver tolerances.
A consequence of this is that solver-tolerance errors in the divergence accumulate over each time step.
\Cref{fig:div} shows this directly:
for CPG(2) at the finest time step $\Delta t$, the $L^2([0,8]; L^2(\Omega))$ norm of $\nabla \cdot \mathbf{u}_h$ reaches around $2 \times 10^{-3}$, three orders of magnitude larger than the corresponding DG(2) value of around $2 \times 10^{-6}$.
A straightforward fix is to modify the temporal quadrature on the incompressibility equation $(\mathrm{div}\,\mathbf{u}_h, q_h) = 0$.
At the discrete level, an $s$-stage quadrature rule enforces $\mathrm{div}\,\mathbf{u}_h = 0$ at every quadrature node up to solver tolerances.
The incompressibility constraint is linear in $\mathbf{u}_h$, so in exact arithmetic any such quadrature rule would be equivalent.
On the discrete level, then, it is advisable to use an $s$-stage quadrature rule over $I_n$ with collocation at $t_{n+1}$, e.g.~a right-Radau rule.
This ensures $\mathrm{div}\,\mathbf{u}_h = 0$ at every time step $t_n$ up to solver tolerances, inhibiting the accumulation of errors.
The resulting stabilized CPG scheme keeps its divergence norm comparable to DG's at all $\Delta t$ (\Cref{fig:div}) and recovers convergence at the finest time step (\Cref{fig:dlerrdt}).

\bda{We can probably make this figure wider, no?}

\begin{figure}
  \begin{minipage}[t]{0.475\textwidth}
    \centering
    \begin{tikzpicture}[scale=0.7]
      \begin{loglogaxis}[
          xlabel={$\Delta t$}, ylabel={Divergence},
          x dir=reverse, ylabel near ticks,
          legend pos=north west,
          legend style={font=\scriptsize, fill opacity=0.85, text opacity=1},
          cycle list name=color list
          ]
          \addplot[seaborngreen, line width = 2, mark=square, mark size = 3]
          table [x=dt,y=L2div, col sep=comma]{data/cylinder/div_dt_dg1.csv};
          \addlegendentry{DG(1)}

          \addplot[seaborngreen, line width = 2, mark=*, mark size = 3]
          table [x=dt,y=L2div, col sep=comma]{data/cylinder/div_dt_dg2.csv};
          \addlegendentry{DG(2)}

          \addplot[seabornred, line width = 2, mark=square, mark size = 3, mark options = {solid}, dashed]
          table [x=dt,y=L2div, col sep=comma]{data/cylinder/div_dt_cpgaux1.csv};
          \addlegendentry{CPG(1)}

          \addplot[seabornred, line width = 2, mark=*, mark size = 3, mark options = {solid}, dashed]
          table [x=dt,y=L2div, col sep=comma]{data/cylinder/div_dt_cpgaux2.csv};
          \addlegendentry{CPG(2)}

          \addplot[seabornorange, line width = 2, mark=square, mark size = 3, mark options = {solid}, dashed]
          table [x=dt,y=L2div, col sep=comma]{data/cylinder/div_dt_cpg_aux_stable1.csv};
          \addlegendentry{stabilized CPG(1)}

          \addplot[seabornorange, line width = 2, mark=*, mark size = 3, mark options = {solid}, dashed]
          table [x=dt,y=L2div, col sep=comma]{data/cylinder/div_dt_cpg_aux_stable2.csv};
          \addlegendentry{stabilized CPG(2)}
      \end{loglogaxis}
    \end{tikzpicture}
    \caption{$L^2([0,8]; L^2(\Omega))$ norm of $\nabla \cdot \mathbf{u}_h$.}
    \label{fig:div}
  \end{minipage}
  \hfill
  \begin{minipage}[t]{0.475\textwidth}
    \centering
    \begin{tikzpicture}[scale=0.7]
      \begin{axis}[xlabel={$\Delta t$}, ylabel={Run time (s)},
          xmode=log, ymode=log, ylabel near ticks,
          legend pos=north east,
          legend style={font=\scriptsize, fill opacity=0.85, text opacity=1},
          cycle list name=color list]
        \addplot[seaborngreen, line width = 2, mark = square, mark size = 3]
        table [x=dt,y=time, col sep=comma]{data/cylinder/stats_dt_dg1.csv};
        \addlegendentry{DG(1)}
        \addplot[seaborngreen, line width = 2, mark = *, mark size = 3]
        table [x=dt,y=time, col sep=comma]{data/cylinder/stats_dt_dg2.csv};
        \addlegendentry{DG(2)}
        \addplot[seabornorange, line width = 2, dashed, mark = square, mark size = 3, mark options = {solid}]
        table [x=dt,y=time, col sep=comma]{data/cylinder/stats_dt_cpg_aux_stable1.csv};
        \addlegendentry{stabilized CPG(1)}
        \addplot[seabornorange, line width = 2, dashed, mark = *, mark size = 3, mark options = {solid}]
        table [x=dt,y=time, col sep=comma]{data/cylinder/stats_dt_cpg_aux_stable2.csv};
        \addlegendentry{stabilized CPG(2)}
      \end{axis}
    \end{tikzpicture}
    \caption{Total run-time versus time step size for DG and stabilized CPG schemes.}
    \label{fig:dlrundt}
  \end{minipage}
\end{figure}

Despite this, DG yields smaller lift and drag errors overall than CPG.
Part of this gap may be because
our reference solution is itself a fine DG simulation, so DG-discretized schemes may benefit from a small inherited bias.
The bulk of the gap, however, is more likely due to DG being better suited than CPG to this problem.

\begin{figure}
  \begin{subfigure}[t]{0.475\textwidth}
    \caption{Newton steps per time step}
    \label{fig:dlnewt}
    \begin{tikzpicture}[scale=0.7]
    \begin{axis}[xlabel={$\Delta t$}, ylabel={Newton iterations},
        xmode = log,
        legend pos=north west, legend style={font=\scriptsize, fill opacity=0.85, text opacity=1},
        cycle list name=color list]
      \addplot[seaborngreen, line width = 2, mark = square, mark size = 3]
      table [x=dt,y=newt, col sep=comma]{data/cylinder/stats_dt_dg1.csv};
      \addlegendentry{DG(1)}
      \addplot[seaborngreen, line width = 2, mark = *, mark size = 3]
      table [x=dt,y=newt, col sep=comma]{data/cylinder/stats_dt_dg2.csv};
      \addlegendentry{DG(2)}
      \addplot[seabornorange, line width = 2, dashed, mark = square, mark size = 3, mark options = {solid}]
      table [x=dt,y=newt, col sep=comma]{data/cylinder/stats_dt_cpg_aux_stable1.csv};
      \addlegendentry{stabilized CPG(1)}
      \addplot[seabornorange, line width = 2, dashed, mark = *, mark size = 3, mark options = {solid}]
      table [x=dt,y=newt, col sep=comma]{data/cylinder/stats_dt_cpg_aux_stable2.csv};
      \addlegendentry{stabilized CPG(2)}
    \end{axis}
    \end{tikzpicture}
  \end{subfigure}
  \begin{subfigure}[t]{0.475\textwidth}
    \caption{Krylov iterations per Newton step}
    \label{fig:dlgmres}
    \begin{tikzpicture}[scale=0.7]
    \begin{axis}[xlabel={$\Delta t$}, ylabel={Krylov iterations},
        xmode = log,
        % legend pos=north west, legend style={font=\scriptsize, fill opacity=0.85, text opacity=1},
        cycle list name=color list]
      \addplot[seaborngreen, line width = 2, mark = square, mark size = 3]
      table [x=dt,y=ksp, col sep=comma]{data/cylinder/stats_dt_dg1.csv};
      % \addlegendentry{DG(1)}
      \addplot[seaborngreen, line width = 2, mark = *, mark size = 3]
      table [x=dt,y=ksp, col sep=comma]{data/cylinder/stats_dt_dg2.csv};
      % \addlegendentry{DG(2)}
      \addplot[seabornorange, line width = 2, dashed, mark = square, mark size = 3, mark options = {solid}]
      table [x=dt,y=ksp, col sep=comma]{data/cylinder/stats_dt_cpg_aux_stable1.csv};
      % \addlegendentry{stabilized CPG(1)}
      \addplot[seabornorange, line width = 2, dashed, mark = *, mark size = 3, mark options = {solid}]
      table [x=dt,y=ksp, col sep=comma]{data/cylinder/stats_dt_cpg_aux_stable2.csv};
      % \addlegendentry{stabilized CPG(2)}
    \end{axis}
    \end{tikzpicture}
  \end{subfigure}
  \caption{Solver performance across time steps for DG and stabilized CPG schemes.}
  \label{fig:dlstatsdt}
\end{figure}
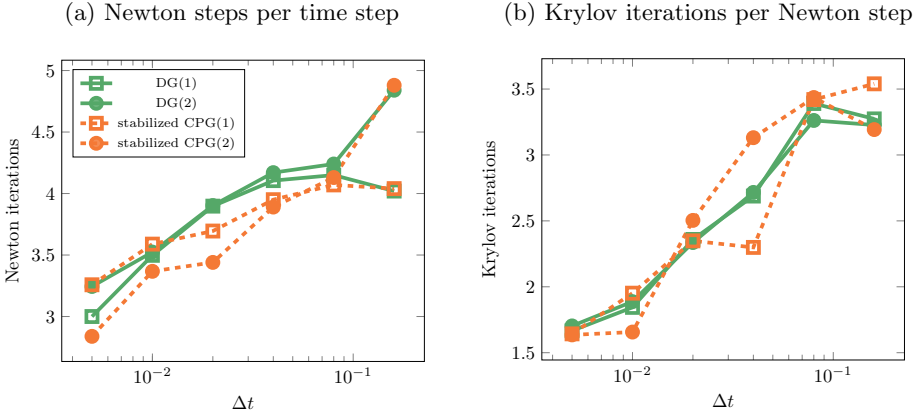

We also report on the performance of our Newton/Krylov/multigrid schemes for the DG and stabilized (with auxiliary variables and under-integrated pressure) CPG schemes.
Our timing experiments were performed on Kodiak, the high-performance computing cluster at Baylor University.  The compute nodes have dual 18-core Intel Xeon Gold 6140 CPUs with 256 GB RAM.  Each run used 16 cores of a single node.
Our solvers are fairly robust with respect to the time step.
\Cref{fig:dlrundt} reports the total time required in the time stepping loop in each of the cases, while \Cref{fig:dlstatsdt} plots the average number of Newton steps per time step and the average number of MG-preconditioned GMRES steps per Newton step.
As larger time steps are taken, we see a slight increase in the number of Newton iterations required per time step and the number of Krylov steps required per Newton step.
However, these are fairly stable so that we see an overall reduction in total run-time with the higher-order method and large time steps.
For example, DG(2) with $\Delta t=0.04$ takes much less time than DG(1) with $\Delta t=0.01$, despite the methods giving comparable overall accuracy.
We also see that the relative timings corresponded to expectations:
the single-stage CPG(1) scheme is the least expensive;
DG(1) and CPG(2) are both more expensive two-stage methods;
the most expensive is the three-stage DG(2) simulation.

\spm{If we're going to give timings here, we should state on what machine the runs were done, and if they were done in parallel.}
\rck{Done.  I put in the information at the beginning of the timings}

\subsection{Compressible flow (Navier--Stokes--Fourier)}
Finally, we implement a novel discretization of compressible flow.
To define the full Navier--Stokes--Fourier system, we introduce the fields in \Cref{tab:nsf_fields}.
The compressible equations are then as follows:
\begin{subequations}\label{eq:nsf}
\begin{align}
    \dot{\rho}  &=  - \nabla \cdot [\rho\mathbf{u}],  \\
    \rho\dot{\mathbf{u}}  &=  - \rho\mathbf{u} \cdot \nabla \mathbf{u} - \nabla p + \frac{2}{\mathrm{Re}}\nabla \cdot \tau,  \\
    \rho\dot{\varepsilon}  &=  -  \nabla \cdot [\varepsilon\mathbf{u}] - p \nabla \cdot \mathbf{u} + \frac{2}{\mathrm{Re}}\tau : \nabla_{\rm sym}\mathbf{u} + \frac{1}{\mathrm{Re}\mathrm{Pr}}\Delta\theta,
\end{align}
\end{subequations}
where the nondimensionalized quantities $\mathrm{Re}$ and $\mathrm{Pr}$ are the Reynolds and Prandtl numbers respectively, and $\nabla_{\rm sym} \mathbf{u} := \frac{1}{2}(\nabla\mathbf{u} + \nabla\mathbf{u}^\top)$ is the symmetric gradient.
\begin{table}
    \centering
    \begin{tabular}{c|l|l}
      Variable & Meaning & Co-domain \\ \hline
      $\rho$ & density & \quad $\mathbb{R}_+$ \\
      $\varepsilon$ & (internal) energy density & \quad $\mathbb{R}_+$ \\
      $p$ & pressure &\quad  $\mathbb{R}$ \\
      $\theta$ & temperature & \quad $\mathbb{R}$ \\
      $\mathbf{u}$ & velocity & \quad $\mathbb{R}^d$ \\
      $\tau$ & stress & \quad $\mathbb{R}^{d \times d}_{\mathrm{sym}}$
    \end{tabular}
    \caption{Fields in the compressible Navier--Stokes--Fourier equations.}
    \label{tab:nsf_fields}
\end{table}

The first four quantities in \Cref{tab:nsf_fields} are coupled by certain equations of state.
In an ideal fluid with non-dimensional heat capacity at constant volume $C_V$,
\begin{subequations}
\begin{equation}
    p = \rho \theta, \qquad
    \varepsilon = C_V p.
\end{equation}
The specific entropy $s$ (such that the generated total entropy is $\int_\Omega \rho s$) can be further related to these quantities;
again, for an ideal fluid
\begin{equation}
    s = C_V \ln \theta - \ln \rho.
\end{equation}
\end{subequations}
The final quantity, the deviatoric stress $\tau$, is a function of the velocity $\mathbf{u}$.
For an isotropic, Newtonian fluid under the Stokes hypothesis, this takes the form
\begin{equation}
  \tau = \nabla_{\rm sym}\mathbf{u} - \frac{1}{3}(\nabla\cdot\mathbf{u})I.
\end{equation}

The variational time discretization considered in~\cite{andrews2025enforcing} proposes introducing certain derived quantities, defined in \Cref{tab:nsf_aux_fields}, to facilitate conservation of mass, momentum, and energy, and the generation of entropy.
For $V_h \subset H^1$, a mixed semidiscretization for the Navier--Stokes--Fourier system \eqref{eq:nsf} may then be defined over these derived fields as:
Find $(\sigma, \boldsymbol{\mu}, \zeta, \tilde{g}, \tilde{\mathbf{u}}, \tilde{\beta}) : \mathbb{R}_+ \to V_h \times V_h^d \times V_h \times V_h \times V_h^d \times V_h$ such that
\begin{subequations}\label{eq:nsf_sd}
\begin{align}
    \int_\Omega 2\sigma \dot{\sigma} v_\rho
        &=  \int_\Omega \tilde{\rho}\tilde{\mathbf{u}} \cdot \nabla v_\rho,  \label{eq:nsf_sd_a} \\
    \int_\Omega \sigma \dot{\boldsymbol{\mu}} \cdot \mathbf{v}_m
        &=  \int_\Omega \frac{1}{2}\tilde{\rho}\tilde{\mathbf{u}} \cdot (\nabla \mathbf{v}_m \cdot \tilde{\mathbf{u}} - \nabla \tilde{\mathbf{u}} \cdot \mathbf{v}_m) - \mathbf{v}_m \cdot \nabla\tilde{p} - \frac{2}{\mathrm{Re}}\tau : \nabla_{\rm sym}\mathbf{v}_m,  \label{eq:nsf_sd_b} \\
    \int_\Omega \varepsilon \dot{\zeta} v_\varepsilon
        &=  \int_\Omega \tilde{\mathbf{u}} \cdot (\tilde{\varepsilon}\nabla v_\varepsilon + \nabla[\tilde{p}v_\varepsilon]) + \frac{2}{\mathrm{Re}}(\tau : \nabla_{\rm sym}\mathbf{u})v_\varepsilon + \frac{1}{\mathrm{Re}\mathrm{Pr}}\theta^2\nabla\tilde{\beta}\cdot\nabla v_\varepsilon,  \label{eq:nsf_sd_c} \\
    \int_\Omega 2\sigma \tilde{g} v_g
        &=  \int_\Omega 2\sigma g v_g,  \label{eq:nsf_sd_d} \\
    \int_\Omega \sigma \tilde{\mathbf{u}} \cdot \mathbf{v}_u
        &=  \int_\Omega \boldsymbol{\mu} \cdot \mathbf{v}_u,  \label{eq:nsf_sd_e} \\
    \int_\Omega \varepsilon \tilde{\beta} v_\beta
        &=  \int_\Omega \varepsilon\beta v_\beta, \label{eq:nsf_sd_f}
\end{align}
\end{subequations}
at all times $t$ and for all $(v_\rho, \mathbf{v}_m, v_\varepsilon, v_g, \mathbf{v}_u, v_\beta) \in V_h \times V_h^d \times V_h \times V_h \times V_h^d \times V_h$.
\begin{table}
    \centering
    \begin{tabular}{c|l}
      Variable & Definition \\ \hline
      $\sigma$ & $\sqrt{\rho}$ \\
      $\zeta$ & $\ln \varepsilon$ \\
      $g$ & $s - \frac{\varepsilon + p}{\rho\theta}$ \\
      $\beta$ & $\frac{1}{\theta}$ \\
      $\boldsymbol{\mu}$ & $\sqrt{\rho} \mathbf{u}$ \\
    \end{tabular}
    \caption{Derived Navier--Stokes--Fourier fields considered in~\cite{andrews2025enforcing}.}
    \label{tab:nsf_aux_fields}
\end{table}

Certain variables, $g$, $\beta$, $\tilde{\rho}$, $\tilde{p}$ and $\tilde{\varepsilon}$, feature in the right-hand side of the discrete system \eqref{eq:nsf_sd}, and require definition in terms of the primal variables $\sigma$, $\zeta$, $\tilde{g}$ and $\tilde{\beta}$.
From $\sigma = \sqrt{\rho}$ and $\zeta = \ln\varepsilon$, one may derive each of the scalar fields in Tables~\ref{tab:nsf_fields} \&~\ref{tab:nsf_aux_fields} via the fluid's equations of state;
in particular $g$ and $\beta$, as featured in \eqref{eq:nsf_sd_d} and \eqref{eq:nsf_sd_f} respectively, can be found for an ideal fluid as
\begin{equation}
    g  =  C_V\zeta - 2(C_V + 1)\ln\sigma - (C_V + 1 + C_V \ln C_V),  \qquad
    \beta  =  C_V \sigma^2\exp(- \zeta).
\end{equation}
One may similarly derive auxiliary approximations (marked with tildes) to the scalar fields from the auxiliary variables $\tilde{g} \approx g = s - \tfrac{\varepsilon + p}{\rho\theta}$ and $\tilde{\beta} \approx \beta = \tfrac{1}{\theta}$;
the fields $\tilde{\rho}$, $\tilde{p}$, $\tilde{\varepsilon}$, as featured in \eqref{eq:nsf_sd_a}, \eqref{eq:nsf_sd_b} and \eqref{eq:nsf_sd_c}, may be found for an ideal fluid as
\begin{equation}
    \tilde{\rho}  =  \tilde{\beta}^{- C_V}\exp(- g - C_V - 1),  \qquad
    \tilde{p}  =  \frac{\tilde{\rho}}{\tilde{\beta}},  \qquad
    \tilde{\varepsilon}  =  C_V\tilde{p}.
\end{equation}

A CPG discretization of \eqref{eq:nsf_sd} was shown in~\cite{andrews2025enforcing} to preserve the conservation of mass $\int_\Omega \sigma^2$, momentum $\int_\Omega \sigma\boldsymbol{\mu}$ and energy $\int_\Omega \tfrac{1}{2}\|\boldsymbol{\mu}\|^2 + \varepsilon$, alongside the generation of entropy $\int_\Omega \rho s$, a property that notably does not hold for a traditional Runge--Kutta discretization.
Since they are not differentiated in time, the variables $\tilde{g}$, $\tilde{\mathbf{u}}$ and $\tilde{\beta}$ are auxiliary variables in the sense outlined in \Cref{sec:aux_vars}.

\begin{figure}[t]
\centering
\begin{tikzpicture}[scale=0.7]
\begin{axis}[
    xlabel={$t$},
    ylabel={$\int_\Omega \rho s$},
    ylabel near ticks,
    legend pos=north west,
    legend style={font=\scriptsize, fill opacity=0.85, text opacity=1},
    grid=major,
    scaled ticks=false,
    x tick label style={/pgf/number format/fixed, /pgf/number format/precision=2}
  ]
  \addplot[line width = 2] table[x=t, y=entropy, col sep=comma]{data/nse/supersonic.dat};
\end{axis}
\end{tikzpicture}
%\caption{Entropy}
%  \end{subfigure}
  \caption{Generation of entropy for the Navier--Stokes--Fourier model.}\label{fig:navier_stokes_fourier}
\end{figure}
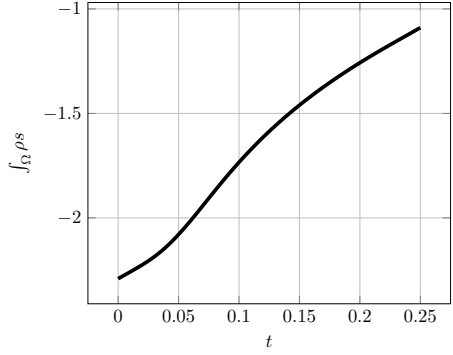

To validate the conservative and dissipative properties of our implementation, we reconsider the test problem from~\cite[Sec.~8.3.2]{andrews2025enforcing} simulating a supersonic perturbation in the velocity field over the periodic unit square $(0, 1)^2$.
Up to interpolation, the initial data is
\begin{subequations}
\begin{align}
    \sigma(0)           &=  1,  \\
    \boldsymbol{\mu}(0) &=  8 \exp(\cos(2\pi x) + \cos(2\pi(y - 0.5)) - 2)\, \mathbf{e}_1,  \\
    \zeta(0)            &=  0,
\end{align}
\end{subequations}
where $\mathbf{e}_1$ is the standard Cartesian basis vector.
We consider an ideal fluid with $C_V = \tfrac{5}{2}$, $\operatorname{Pr} = 0.71$ and $\operatorname{Re} = 128$.
The spatial discretization uses a $256 \times 256$ quadrilateral mesh with continuous bilinear elements in space.
We use a 1-stage CPG discretization in time (with auxiliary variables correspondingly in the degree-0 DG space) over a uniform time step $\Delta t = 1/2048$.
\Cref{fig:navier_stokes_fourier} shows the generation of entropy, where we note that we achieve (as expected) conservation of mass, momentum and energy up to solver tolerances and roundoff error.

\section{Conclusions and future work}\label{sec:conclusions}

Galerkin-in-time discretizations have recently attracted substantial interest for their structure-preserving properties, relative to classical implicit \RK{} methods.
We have extended the \Irksome{} package to implement them within the same interface used for \RK{} discretizations.
The resulting package allows easy selection of both CPG and DG discretizations, flexible customization of time quadrature, auxiliary variables for CPG discretizations of differential--algebraic equations, and stage-coupled algebraic solvers including monolithic multigrid.
The numerical examples provided demonstrate the package is simple to use within the standard \Irksome{} interface across a range of problems, and confirm the expected structure-preserving properties are realized in practice.
While we hope the software will support further study and development of structure-preserving time stepping schemes, we note that open questions remain, notably around the development of efficient linear and nonlinear solvers tailored to these schemes.

\bibliographystyle{siamplain}
\bibliography{references}

@misc{BarattaEtal2023,
  title     = {{DOLFINx}: the next generation {FEniCS} problem solving environment},
  author    = {Baratta, Igor A. and Dean, Joseph P. and Dokken, J{\o}rgen S. and Habera, Michal and Hale, Jack S. and Richardson, Chris N. and Rognes, Marie E. and Scroggs, Matthew W. and Sime, Nathan and Wells, Garth N.},
  doi       = {10.5281/zenodo.10447666},
  year      = {2023},
  howpublished = {preprint}
}

@article{shen2018scalar,
  title={The scalar auxiliary variable ({SAV}) approach for gradient flows},
  author={Shen, Jie and Xu, Jie and Yang, Jiang},
  journal={Journal of Computational Physics},
  volume={353},
  pages={407--416},
  year={2018},
  publisher={Elsevier},
  doi = {10.1016/j.jcp.2017.10.021},
}

@article{estep2002dynamical,
  title={The dynamical behavior of the discontinuous {G}alerkin method and related difference schemes},
  author={Estep, Donald and Stuart, Andrew},
  journal = {Mathematics of Computation},
  volume  = {71},
  number  = {239},
  pages   = {1075--1103},
  year    = {2002},
  doi     = {10.1090/S0025-5718-01-01364-3}
}

@article{bao2004computing,
  title={Computing the ground state solution of {Bose--Einstein} condensates by a normalized gradient flow},
  author={Bao, Weizhu and Du, Qiang},
  journal={SIAM Journal on Scientific Computing},
  volume={25},
  number={5},
  pages={1674--1697},
  year={2004},
  publisher={SIAM},
  doi={10.1137/S1064827503422956},
}

@misc{abhyankar2018petsc,
  title         = {{PETSc/TS: A modern scalable ODE/DAE solver library}},
  author        = {Abhyankar, Shrirang and Brown, Jed and Constantinescu, Emil M. and Ghosh, Debojyoti and Smith, Barry F. and Zhang, Hong},
  year          = {2018},
  eprint        = {1806.01437},
  archiveprefix = {arXiv},
  primaryclass  = {math.NA},
}

@article{hindmarsh2005sundials,
  title     = {{SUNDIALS}: Suite of nonlinear and differential/algebraic equation solvers},
  author    = {Hindmarsh, Alan C. and Brown, Peter N. and Grant, Keith E. and Lee, Steven L. and Serban, Radu and Shumaker, Dan E. and Woodward, Carol S.},
  journal   = {ACM Transactions on Mathematical Software},
  volume    = {31},
  number    = {3},
  pages     = {363--396},
  year      = {2005},
  doi       = {10.1145/1089014.1089020},
  publisher = {ACM New York, NY, USA}
}

@article{abu2022monolithic,
  title   = {Monolithic multigrid for implicit {R}unge--{K}utta discretizations of incompressible fluid flow},
  author  = {Abu-Labdeh, Razan and MacLachlan, Scott and Farrell, Patrick E.},
  journal = {Journal of Computational Physics},
  doi     = {10.1016/j.jcp.2023.111961},
  year    = {2023},
  volume  = {478},
  pages   = {111961}
}

@article{alfeld-sorokina16,
  author  = {Alfeld, Peter and Sorokina, Tatyana},
  title   = {Linear differential operators on bivariate spline spaces and spline vector fields},
  journal = {BIT Numerical Mathematics},
  volume  = {56},
  year    = {2016},
  doi     = {10.1007/s10543-015-0557-x},
  pages   = {15--32}
}

@article{Alnaes:2014,
  author        = {Aln\ae{}s, Martin S. and Logg, Anders and \O{}lgaard, Kristian
                   B. and Rognes, Marie E. and Wells, Garth
                   N.},
  title         = {{Unified Form Language: a domain-specific language
                   for weak formulations of partial differential
                   equations}},
  journal       = {ACM Transactions on Mathematical Software},
  year          = 2014,
  volume        = 40,
  number        = 2,
  pages         = {9:1--9:37},
  doi           = {10.1145/2566630},
}

@misc{Andrews_Farrell_2025b,
  title     = {Conservative and dissipative discretisations of multi-conservative {ODEs} and {GENERIC} systems},
  doi       = {10.48550/arXiv.2511.23266},
  publisher = {arXiv},
  author    = {Andrews, B. D. and Farrell, P. E.},
  month     = nov,
  year      = {2025}
}

@article{andrews2025enforcing,
  title   = {Enforcing conservation laws and dissipation inequalities numerically via auxiliary variables},
  author  = {Andrews, Boris D. and Farrell, Patrick E.},
  journal = {SIAM Journal on Scientific Computing},
  volume = {47},
  number = {6},
  pages = {A3516--A3535},
  year = {2025},
  doi = {10.1137/25M1756673},
}

@article{Betsch_Steinmann_2000a,
  title   = {Inherently energy conserving time finite elements for classical mechanics},
  volume  = {160},
  issn    = {0021-9991},
  doi     = {10.1006/jcph.2000.6427},
  number  = {1},
  journal = {Journal of Computational Physics},
  author  = {Betsch, P. and Steinmann, P.},
  month   = may,
  year    = {2000},
  pages   = {88--116}
}

@article{Betsch_Steinmann_2000b,
  title   = {Conservation properties of a time {FE} method—part {I}: time-stepping schemes for {N}-body problems},
  volume  = {49},
  issn    = {1097-0207},
  doi     = {10.1002/1097-0207(20001020)49:5<599::AID-NME960>3.0.CO;2-9},
  number  = {5},
  journal = {International Journal for Numerical Methods in Engineering},
  author  = {Betsch, P. and Steinmann, P.},
  year    = {2000},
  pages   = {599--638}
}

@article{betteridge2024ngspetsc,
  title   = {{ngsPETSc: A coupling between NETGEN/NGSolve and PETSc}},
  author  = {Betteridge, Jack and Farrell, Patrick E and Hochsteger, Matthias and Lackner, Christopher and Sch{\"o}berl, Joachim and Zampini, Stefano and Zerbinati, Umberto},
  journal = {Journal of Open Source Software},
  volume  = {9},
  number  = {104},
  pages   = {7359},
  year    = {2024},
  doi = {10.21105/joss.07359},
}

@misc{brubeck2025fiat,
  title         = {{FIAT}: enabling classical and modern macroelements},
  author        = {Brubeck, Pablo D. and Kirby, Robert C.},
  year          = {2025},
  eprint        = {2501.14599},
  archiveprefix = {arXiv},
  primaryclass  = {math.NA},
}

@article{Celledoni_et_al_2009,
  title     = {Energy-preserving {Runge}--{Kutta} methods},
  author    = {Celledoni, E. and McLachlan, R. I. and McLaren, D. I. and Owren, B. and Quispel, G. R. W. and Wright, W. M.},
  journal   = {ESAIM: Mathematical Modelling and Numerical Analysis},
  volume    = {43},
  number    = {4},
  pages     = {645--649},
  year      = {2009},
  publisher = {EDP Sciences},
  doi       = {10.1051/m2an/2009020}
}

@article{chen2025unconditionally,
  title     = {Unconditionally energy stable {IEQ-FEM}s for the {C}ahn--{H}illiard equation and {A}llen--{C}ahn equation},
  author    = {Chen, Yaoyao and Liu, Hailiang and Yi, Nianyu and Yin, Peimeng},
  journal   = {Numerical Algorithms},
  volume    = {99},
  number    = {3},
  pages     = {1161--1202},
  year      = {2025},
  publisher = {Springer},
  doi       = {10.1007/s11075-024-01910-z}
}

@article{Cohen_Hairer_2011,
  title    = {Linear energy-preserving integrators for {Poisson} systems},
  volume   = {51},
  issn     = {1572-9125},
  doi      = {10.1007/s10543-011-0310-z},
  number   = {1},
  journal  = {BIT Numerical Mathematics},
  author   = {Cohen, D. and Hairer, E.},
  month    = mar,
  year     = {2011},
  pages    = {91--101}
}

@article{Egger_Habrich_Shashkov_2021,
  title    = {On the energy stable approximation of {Hamiltonian} and gradient systems},
  volume   = {21},
  issn     = {1609-9389},
  doi      = {10.1515/cmam-2020-0025},
  number   = {2},
  journal  = {Computational Methods in Applied Mathematics},
  author   = {Egger, H. and Habrich, O. and Shashkov, V.},
  month    = apr,
  year     = {2021},
  pages    = {335--349}
}

@book{Ern_Guermond_2021c,
  address    = {Cham, Switzerland},
  series     = {Texts in Applied Mathematics},
  title      = {Finite Elements {III}: First-Order and Time-Dependent {PDEs}},
  volume     = {74},
  isbn       = {978-3-030-57347-8 978-3-030-57348-5},
  shorttitle = {Finite Elements {III}},
  publisher  = {Springer International Publishing},
  author     = {Ern, A. and Guermond, J.-L.},
  year       = {2021},
  doi        = {10.1007/978-3-030-57348-5}
}

@article{farrell2021irksome,
  title   = {{Irksome: Automating Runge--Kutta time-stepping for finite element methods}},
  author  = {Farrell, Patrick E. and Kirby, Robert C. and Marchena-Menendez, Jorge},
  volume  = {47},
  number  = {4},
  pages   = {1--26},
  year    = {2021},
  doi     = {10.1145/3466168},
  journal = {ACM Transactions on Mathematical Software}
}

@article{farrell2021pcpatch,
  title     = {{PCPATCH}: software for the topological construction of multigrid relaxation methods},
  author    = {Farrell, Patrick E. and Knepley, Matthew G. and Mitchell, Lawrence and Wechsung, Florian},
  journal   = {ACM Transactions on Mathematical Software},
  volume    = {47},
  number    = {3},
  pages     = {1--22},
  year      = {2021},
  publisher = {ACM New York, NY, USA},
  doi = {10.1145/3445791},
}

@article{French_Schaeffer_1990,
  title   = {Continuous finite element methods which preserve energy properties for nonlinear problems},
  volume  = {39},
  issn    = {0096-3003},
  doi     = {10.1016/S0096-3003(20)80006-X},
  number  = {3},
  journal = {Applied Mathematics and Computation},
  author  = {French, D. A. and Schaeffer, J. W.},
  month   = oct,
  year    = {1990},
  pages   = {271--295}
}

@article{Giesselmann_Karsai_Tscherpel_2025,
  title   = {Energy-consistent {Petrov}--{Galerkin} time discretization of port-{Hamiltonian} systems},
  volume  = {11},
  issn    = {2426-8399},
  doi     = {10.5802/smai-jcm.127},
  journal = {SMAI Journal of Computational Mathematics},
  author  = {Giesselmann, J. and Karsai, A. and Tscherpel, T.},
  year    = {2025},
  pages   = {335--367}
}

@article{guzman2018inf,
  title     = {Inf-sup stable finite elements on barycentric refinements producing divergence--free approximations in arbitrary dimensions},
  author    = {Guzm{\'a}n, Johnny and Neilan, Michael},
  journal   = {SIAM Journal on Numerical Analysis},
  volume    = {56},
  number    = {5},
  pages     = {2826--2844},
  year      = {2018},
  publisher = {SIAM},
  doi = {10.1137/17M1153467},
}

@article{Hairer_Lubich_2014,
  title    = {Energy-diminishing integration of gradient systems},
  volume   = {34},
  issn     = {0272-4979},
  doi      = {10.1093/imanum/drt031},
  number   = {2},
  journal  = {IMA Journal of Numerical Analysis},
  author   = {Hairer, E. and Lubich, C.},
  month    = apr,
  year     = {2014},
  pages    = {452--461}
}

@book{hairer2006geometric,
  title     = {Geometric numerical integration: structure-preserving algorithms for ordinary differential equations},
  author    = {Hairer, Ernst and Lubich, Christian and Wanner, Gerhard},
  volume    = {31},
  year      = {2006},
  publisher = {Springer Science \& Business Media}
}

@techreport{huynh2011collocation,
  author      = {Huynh, H. T.},
  title       = {Collocation and {G}alerkin time-stepping methods},
  institution = {NASA Glenn Research Center},
  number      = {NASA/TM-2011-216340},
  year        = {2011},
  url         = {https://ntrs.nasa.gov/citations/20110014969}
}

@article{john2004reference,
  title     = {Reference values for drag and lift of a two-dimensional time-dependent flow around a cylinder},
  author    = {John, Volker},
  journal   = {International Journal for Numerical Methods in Fluids},
  volume    = {44},
  number    = {7},
  pages     = {777--788},
  year      = {2004},
  publisher = {Wiley Online Library},
  doi = {10.1002/fld.679},
}

@book{kevrekidis2015,
  doi       = {10.1137/1.9781611973945},
  year      = 2015,
  publisher = {SIAM},
  author    = {Panayotis G. Kevrekidis and Dimitri J. Frantzeskakis and Ricardo Carretero-Gonz{\'{a}}lez},
  title     = {The Defocusing Nonlinear Schrödinger Equation}
}

@article{Kirby:2004,
  author  = {Kirby, Robert C.},
  title   = {{Algorithm 839: FIAT, a new paradigm for computing
             finite element basis functions}},
  journal = {ACM Transactions on Mathematical Software},
  year    = 2004,
  volume  = 30,
  number  = 4,
  pages   = {502--516},
  doi     = {10.1145/1039813.1039820}
}

@article{kirby2018solver,
  author    = {Kirby, Robert C. and Mitchell, Lawrence},
  title     = {Solver composition across the {PDE}/linear algebra
               barrier},
  journal   = {SIAM Journal on Scientific Computing},
  year      = 2018,
  volume    = 40,
  number    = 1,
  pages     = {C76--C98},
  publisher = {SIAM},
  doi       = {10.1137/17M1133208}
}

@misc{kirby2025automated,
  title         = {Automated {R}unge--{K}utta--{N}ystr\"{o}m time stepping for finite element methods in {I}rksome},
  author        = {Kirby, Robert C. and MacLachlan, Scott P. and Brubeck, Pablo D.},
  year          = {2025},
  eprint        = {2508.20255},
  archiveprefix = {arXiv},
  primaryclass  = {math.NA},
}

@article{kirby2025extending,
  title     = {Extending {I}rksome: improvements in automated {Runge--Kutta} time stepping for finite element methods},
  author    = {Kirby, Robert C. and MacLachlan, Scott P.},
  journal   = {ACM Transactions on Mathematical Software},
  volume    = {51},
  number    = {3},
  pages     = {1--27},
  year      = {2025},
  doi       = {10.1145/3759245},
  publisher = {ACM New York, NY}
}

@article{McLachlan_Quispel_Robidoux_1999,
  title    = {Geometric integration using discrete gradients},
  volume   = {357},
  doi      = {10.1098/rsta.1999.0363},
  number   = {1754},
  journal  = {Philosophical Transactions of the Royal Society of London. Series A: Mathematical, Physical and Engineering Sciences},
  author   = {McLachlan, R. I. and Quispel, G. R. W. and Robidoux, N.},
  month    = apr,
  year     = {1999},
  pages    = {1021--1045}
}

@article{mmg,
  title     = {On the convergence of monolithic multigrid for implicit {R}unge--{K}utta time stepping of finite element problems},
  author    = {Kirby, Robert C},
  journal   = {SIAM Journal on Scientific Computing},
  volume    = {46},
  number    = {5},
  pages     = {S22--S45},
  year      = {2024},
  publisher = {SIAM},
  doi = {10.1137/23M1569344},
}

@article{Pavarino:1993,
  author  = {Luca F. Pavarino},
  title   = {{Additive Schwarz methods for the $p$-version finite
             element method}},
  journal = {Numerische Mathematik},
  year    = 1993,
  volume  = 66,
  number  = 1,
  pages   = {493--515},
  issn    = {0945-3245},
  doi     = {10.1007/BF01385709}
}

@manual{FiredrakeUserManual,
  title        = {Firedrake User Manual},
  author       = {David A. Ham and Paul H. J. Kelly and Lawrence Mitchell and Colin J. Cotter and Robert C. Kirby and Koki Sagiyama and Nacime Bouziani and Sophia Vorderwuelbecke and Thomas J. Gregory and Jack Betteridge and Daniel R. Shapero and Reuben W. Nixon-Hill and Connor J. Ward and Patrick E. Farrell and Pablo D. Brubeck and India Marsden and Thomas H. Gibson and Miklós Homolya and Tianjiao Sun and Andrew T. T. McRae and Fabio Luporini and Alastair Gregory and Michael Lange and Simon W. Funke and Florian Rathgeber and Gheorghe-Teodor Bercea and Graham R. Markall},
  organization = {Imperial College London and University of Oxford and Baylor University and University of Washington},
  edition      = {First edition},
  year         = {2023},
  month        = {5},
  doi          = {10.25561/104839},
}

@incollection{schafer1996benchmark,
  title     = {Benchmark computations of laminar flow around a cylinder},
  author    = {Sch{\"a}fer, Michael and Turek, Stefan and Durst, Franz and Krause, Egon and Rannacher, Rolf},
  booktitle = {Flow simulation with high-performance computers II},
  pages     = {547--566},
  year      = {1996},
  publisher = {Springer},
  doi = {10.1007/978-3-322-89849-4_39},
}

@article{Schoeberl:2008,
  author  = {Joachim Sch\"oberl and Jens M. Melenk and Clemens
             Pechstein and Sabine Zaglmayr},
  title   = {{Additive Schwarz preconditioning for $p$-version
             triangular and tetrahedral finite elements}},
  journal = {IMA Journal of Numerical Analysis},
  year    = 2008,
  volume  = 28,
  pages   = {1--24},
  doi     = {10.1093/imanum/drl046}
}

@techreport{scott1984conforming,
  title       = {Conforming Finite Element Methods for Incompressible and Nearly Incompressible Continua},
  author      = {Scott, L.~Ridgway and Vogelius, Michael},
  year        = {1984},
  institution = {Maryland Univ.~College Park Inst. For Physical Science And Technology}
}

@article{vanka1986block,
  title     = {Block-implicit multigrid solution of {N}avier--{S}tokes equations in primitive variables},
  author    = {Vanka, S.~Pratap},
  journal   = {Journal of Computational Physics},
  volume    = {65},
  number    = {1},
  pages     = {138--158},
  year      = {1986},
  publisher = {Elsevier},
  doi = {10.1016/0021-9991(86)90008-2},
}

@article{vanlent2005,
  doi     = {10.1137/030601144},
  year    = 2005,
  volume  = {27},
  number  = {1},
  pages   = {67--92},
  author  = {Van Lent, J. and Vandewalle, S.},
  title   = {Multigrid Methods for Implicit {Runge--Kutta} and Boundary Value Method Discretizations of Parabolic {PDEs}},
  journal = {SIAM Journal on Scientific Computing}
}

@misc{zenodo/Zenodo-20260619.7,
 key   = {zenodo/Zenodo-20260619.7},
 title = {{Software used in `Automated Galerkin time stepping in Irksome'}},
 year  = {2026},
 month = {jun},
 doi   = {10.5281/zenodo.20784213},
 url   = {https://doi.org/10.5281/zenodo.20784213},
}

@article{allen1972,
  title   = {Ground state structures in ordered binary alloys with second neighbor interactions},
  volume  = {20},
  doi     = {10.1016/0001-6160(72)90037-5},
  number  = {3},
  journal = {Acta Metallurgica},
  author  = {Allen, S. M. and Cahn, J. W.},
  year    = {1972},
  pages   = {423--433}
}

@article{allen1979,
  title   = {A microscopic theory for antiphase boundary motion and its application to antiphase domain coarsening},
  journal = {Acta Metallurgica},
  volume  = {27},
  number  = {6},
  pages   = {1085--1095},
  year    = {1979},
  doi     = {10.1016/0001-6160(79)90196-2},
  author  = {Allen, S. M. and Cahn, J. W.}
}

@article{marsden2001,
	Author = {Marsden, J. E. and West, M.},
	Journal = {Acta Numerica},
	Pages = {357--514},
	Publisher = {Cambridge Univ Press},
	Title = {Discrete mechanics and variational integrators},
	Volume = {10},
	Year = {2001},
    doi = {10.1017/S096249290100006X},
}

@misc{rafiei2025improvingpatchselectionmonolithic,
      title={Improving patch selection for monolithic multigrid solvers for high-order {T}aylor--{H}ood discretizations},
      author={Amin Rafiei and Scott MacLachlan},
      year={2025},
      eprint        = {2502.01130},
      archiveprefix = {arXiv},
      primaryclass  = {math.NA},
}

@article{BASTIAN202175,
title = {The {DUNE} framework: basic concepts and recent developments},
journal = {Computers \& Mathematics with Applications},
volume = {81},
pages = {75--112},
year = {2021},
issn = {0898-1221},
doi = {10.1016/j.camwa.2020.06.007},
author = {Peter Bastian and Markus Blatt and Andreas Dedner and Nils-Arne Dreier and Christian Engwer and Ren\'{e} Fritze and Carsten Gr\"{u}ser and Christoph Gr\"{u}ninger and Dominic Kempf and Robert Kl\"{o}fkorn and Mario Ohlberger and Oliver Sander},
}

@article{SMAI-JCM_2023__9__61_0,
     author = {Versbach, Lea Miko and Linders, Viktor and Kl\"ofkorn, Robert and Birken, Philipp},
     title = {Theoretical and {Practical} {Aspects} of {Space-Time} {DG-SEM} {Implementations}},
     journal = {SMAI Journal of Computational Mathematics},
     pages = {61--93},
     year = {2023},
     publisher = {Soci\'et\'e de Math\'ematiques Appliqu\'ees et Industrielles},
     volume = {9},
     doi = {10.5802/smai-jcm.95},
     language = {en},
}

\end{document}